\documentclass[12pt]{amsart}

\usepackage{a4wide}

\usepackage{amssymb}
\usepackage{epsfig}
\usepackage{subfigure}

\usepackage{enumerate}

\usepackage{array}
\newtheorem{thm}{Theorem}[section]
\newtheorem{dfn}[thm]{Definition}
\newtheorem{cor}[thm]{Corollary}
\newtheorem{prop}[thm]{Proposition}
\newtheorem{lem}[thm]{Lemma}
\theoremstyle{remark} 
\newtheorem{rmk}[thm]{Remark}
\newtheorem{assume}[thm]{\textsc{Assumption}}

\def\cqfd{\mbox{}\nolinebreak\hfill$\Box$\medbreak\par}
\newenvironment{pf}{\noindent\textbf{Proof:}}{\cqfd}

\def\rh#1{{{\bf H}^{#1}_{\R}}}
\def\ch#1{{{\bf H}^{#1}_{\C}}}
\def\chb#1{{{\overline{\bf H}}^{#1}_{\C}}}
\newcommand{\Z}{\mathbb{Z}}
\newcommand{\N}{\mathbb{N}}
\newcommand{\Q}{\mathbb{Q}}
\newcommand{\C}{\mathbb{C}}
\newcommand{\R}{\mathbb{R}}
\newcommand{\B}{\mathcal{B}}
\newcommand{\Ss}{\mathcal{S}}

\renewcommand{\l}{\lambda}
\newcommand{\lc}{\overline{\lambda}}
\renewcommand{\k}{k}
\newcommand{\re}{\mathfrak{Re}}

\newcommand{\presGroup}[2]{\langle\, #1\, |\, #2 \, \rangle}

\title[Spherical CR uniformization]{A 1-parameter family of spherical
  CR uniformizations of the figure eight knot complement}

\author{Martin Deraux}

\date{March 18, 2016}

\begin{document}

\begin{abstract}
  We describe a simple fundamental domain for the holonomy group of
  the boundary unipotent spherical CR uniformization of the figure
  eight knot complement, and deduce that small deformations of that
  holonomy group (such that the boundary holonomy remains parabolic)
  also give a uniformization of the figure eight knot complement.
  Finally, we construct an explicit 1-parameter family of deformations
  of the boundary unipotent holonomy group such that the boundary holonomy
  is twist-parabolic. For small values of the twist of these parabolic
  elements, this produces a 1-parameter family of pairwise
  non-conjugate spherical CR uniformizations of the figure eight knot
  complement.
\end{abstract}

\maketitle

\section{Introduction}

The existence of a complete hyperbolic structure on a 3-manifold has
important topological consequences. For instance, this gives a
definition of the volume of a knot (when a knot admits a complete
hyperbolic structure, that structure is unique by Mostow rigidity, so
the volume of that metric is a well-defined invariant).

In this paper, we focus on another kind of geometric structures on
3-manifolds, namely structures modeled on the boundary of a symmetric
space $X$ of negative curvature (transition maps are required to be
locally given by isometries of $X$). The visual boundary
$\partial_\infty X$ is then a 3-dimensional sphere if
$X=H^4_\R$ or $H^2_\C$.

The first case gives rise to the theory of flat conformal structures,
and the second one to the theory spherical CR structures. In the first
case, one considers the unit ball model of $H^4_\R$, so the visual
boundary is $S^3\subset \R^4$, and the group of isometries of $H^4_\R$
acts as M\"obius transformations (i.e. transformations that map
spheres into spheres, of possibly infinite radius). Alternatively, one
can use stereographic projection and think of $S^3$ as
$\R^3\cup\{\infty\}$; this would also correspond to using the upper
half plane model for $H^3_\R$.

In the second case, using the ball model $\mathbb{B}^2\subset\C^2$ one
can identify $\partial_\infty H^2_\C$ with the unit sphere
$S^3\subset\C^2$. The action on the boundary is best understood in
stereographic projection, and identifying
$S^3\setminus\{p_\infty\}\simeq \R^3\simeq \C\times \R$ with the
Heisenberg group. Isometries of $H^2_C$ fixing $p_\infty$ then acts as
automorphisms of the Heisenberg group.  Of course the Heisenberg group
acting on itself by left translations gives many automorphisms (which
correspond to the action of unipotent matrices in $U(2,1)$), and one
gets the full automorphism group by adjoining a rotation in
$\C\times\R$ around the $\R$ factor, and scaling of the form
$(z,t)\mapsto (\lambda z,\lambda^2t)$ (which corresponds to a
loxodromic isometry), see section~\ref{sec:siegel}.

Even though a lot of partial results have been obtained
(see~\cite{kamishimatsuboi},~\cite{goldmantorus} for instance), the
classification of 3-manifolds that admit a spherical CR structure is
far from understood. When a manifold admits a spherical CR structure,
the moduli space of such structures is also quite mysterious.

In this paper, we will be interested in a special kind of spherical CR
structures, namely spherical CR \emph{uniformizations} (in the
literature, these are sometimes called complete spherical CR
structures). These are characterized by the fact that the developing
map of the structure is a diffeomorphism onto its image, which is an
open set in $S^3$. In that case, the holonomy group is a discrete
subgroup $\Gamma\subset PU(2,1)$, and the image of the developing map
is the domain of discontinuity $\Omega_\Gamma$ of $\Gamma$ (i.e. the
largest open set where the action is proper). The quotient
$\Gamma\setminus\Omega_\Gamma$ is called the \emph{manifold at
  infinity} of $\Gamma$.

The classification of 3-manifolds that admit a spherical CR
uniformization is also an open problem. Recall that $H^2_\C$ is a
homogeneous under the action of $PU(2,1)$, and the isotropy group of a
point is isomorphic to $U(2)$. In particular, finite subgroups of
$U(2)$ such that nontrivial elements fix only the origin (in other
words the groups should not contain any complex reflection) yield
spherical CR uniformizable 3-manifolds with finite fundamental group.

In a similar vein, quotients of the Heisenberg group yield Nil
manifolds that trivially admit a spherical CR uniformization, such
that the holonomy group has a global fixed point, which is now in
$\partial_\infty H^2_\C$ instead of $H^2_\C$.

It is also natural to consider stabilizers of totally geodesic
subspaces in $H^2_\C$, namely copies of $H^2_\R$ or $H^1_\C$. In that
setting, Fuchsian groups (i.e. discrete subgroups of $SO(2,1)$ or
$SU(1,1)$, seen as subgroups of $SU(2,1)$) produce as their manifold
at infinity a circle bundle over a surface (or more generally over a
2-orbifold). This class is more interesting than the previous one,
because it is known that the corresponding groups often admit
deformations (but not always, see~\cite{toledoJDG}).  We will
summarize the results in this well developed line of research by
saying simply that many Seifert 3-manifolds admit spherical CR
uniformizations
(see~\cite{gkl},~\cite{agg},~\cite{parkerplatis},~\cite{will} and
others).

The class of \emph{hyperbolic} manifolds that admit a spherical CR
uniformization is also far from being understood.  In a number of beautiful
results that appeared in the last decade, Schwartz discovered that
many hyperbolic manifolds admit spherical CR uniformizations,
see~\cite{richWLC},~\cite{rhochi} and~\cite{richBook}. His starting
point was to consider representations of triangle groups into
$PU(2,1)$, see~\cite{schwartzICM}, and to determine the manifold at
infinity of well chosen such representations.

More recently, the figure eight knot complement was shown to admit a
spherical CR uniformization~\cite{derauxfalbel} by following a
somewhat different strategy, namely it was found as a byproduct of
Falbel's program for finding representations of fundamental groups of
triangulated 3-manifold into $PU(2,1)$ (see~\cite{falbelFigure8}), or
in $PGL(3,\C)$ (see~\cite{bfg}).

Falbel's construction turned out to produce lots of representations,
and in fact so many that the geometric properties of the resulting
representations are in general difficult to analyze. In order to make
the list more tractable (and also for other reasons related to the
study of Bloch groups), the search is often restricted to
representations such that peripheral subgroups are mapped to unipotent
matrices (matrices with 1 as their only eigenvalue).
The boundary unipotent representations for non-compact 3-manifolds
with low complexity (i.e those that can be built by gluing up to three
ideal tetrahedra) are listed in~\cite{fkr}, and the geometry of some
of these representations are analyzed in~\cite{derauxfalbel}
and~\cite{derauxcensusfkr}. It turns out very few representations in
that list are discrete.

It is quite clear however that the unipotent restriction is somewhat
artificial. Part of the point of the present paper is to show that, at
least in some cases, there are many boundary parabolic representations
that are not unipotent, and that these representations carry just as
much interesting geometric information about the 3-manifold.

Let $M$ denote the figure eight knot complement. The main goal of this
paper is to show that $M$ admits a 1-parameter family of pairwise non
conjugate spherical CR uniformizations.  

We will build on the fact that $M$ admits a unique spherical CR
uniformization with unipotent boundary holonomy, as was shown
in~\cite{derauxfalbel}. For future reference, we will refer to that
structure simply as \emph{the} boundary unipotent uniformization of
$M$ (see the precise uniqueness statement in~\cite{derauxfalbel}), and
we denote the corresponding holonomy representation by $\rho$.  In
view of Schwartz's spherical CR Dehn surgery theorem~\cite{richBook},
one expects that small deformations of the boundary unipotent holonomy
representation should still be discrete, and they should have a
manifold at infinity given by \emph{some} Dehn filling of the figure
eight knot complement.

In order to turn this into a proof, one could try and prove that the
boundary unipotent representation satisfies the hypotheses of
Schwartz's theorem, i.e. that its image is a horotube group (without
exceptional parabolic elements), and that its limit set is porous. If
that works, then it is enough to show that the group admits
deformations, and to study the type of the deformed unipotent element;
Schwartz's surgery formula shows in particular that (under some
technical assumptions), if there are deformations where the unipotent
peripheral holonomy stays parabolic, then the manifold at infinity
should not change at all in small deformations.

Although a few examples of non-compact hyperbolic manifolds
are known to admit spherical CR uniformizations
(see~\cite{richWLC},~\cite{rhochi},~\cite{derauxfalbel}), the
deformation theory of the holonomy representations of these examples
is still quite mysterious. In particular, there are only two examples
where non-trivial deformations are known to exist such that peripheral
elements map to parabolic elements. These two examples are the figure
eight knot complement and the Whitehead link complement. The results
announced by Parker and Will, see~\cite{parkerwill} say that there are
at least two different spherical CR uniformizations of the Whitehead
link complement, and that there is a 1-parameter family of
representations interpolating between their holonomy representations.

The main result of our paper gives an explicit construction of
twist-parabolic deformations.
\begin{thm}\label{thm:param}
  There is a continuous 1-parameter family of irreducible
  representations $\rho_t:\pi_1(M)\rightarrow PU(2,1)$, such that for
  each $t$, $\rho_t$ maps peripheral subgroups of $M$ onto a cyclic
  group generated by a single parabolic element with eigenvalues
  $e^{it},e^{it},e^{-2it}$.
\end{thm}
Given the eigenvalue condition, it should be clear that the
representations $\rho_t$ are pairwise non conjugate. We will choose
$\rho_t$ so that $\rho_0$ is the holonomy of the boundary unipotent
spherical CR uniformization.

Note that the existence of such parabolic deformations was independtly
discovered by Pierre-Vincent Koseleff, using a variant of the method
devised by Falbel to parametrize boundary unipotent representations of
3-manifolds, see~\cite{falbelFigure8},~\cite{bfg} and~\cite{fkr} for
instance. An alternative parametrization of this family can also be
obtained from the description of the full character variety,
see~\cite{fgkrt}, see also~\cite{hmp}.

We will use a more na\"ive construction, which is closer in spirit to
Riley's parametrization of the character variety of the figure eight
knot group (or more generally 2-bridge knot groups) into $PSL_2(\C)$,
see~\cite{riley2bridge}.

Our main result is the following.
\begin{thm}\label{thm:main}
 There exists a $\delta>0$ such that for $|t|<\delta$, $\rho_t$ is the
 holonomy of a spherical CR uniformization of the figure eight knot
 complement.
\end{thm}
In order to show this, we will study the Ford domain for the image of
$\rho_0$, and we will show that it is generic enough for its
combinatorics to be preserved under small deformations of $\rho_0$.
Note that this argument turns out to fail for the Ford domain of the
holonomy of the spherical CR uniformization of the Whitehead link
complement announced by Parker and Will in~\cite{parkerwill}.  Indeed,
their Ford domain has the same local combinatorial structure as the
Dirichlet domain described in~\cite{derauxfalbel}, and in particular
it has lots of tangent spinal spheres.

It will be clear to the reader familiar with the notion of
horotubes~\cite{richBook} that the Ford domain exhibits an explicit
horotube structure for the group, but since our construction of
horotubes is actually very close to proving Theorem~\ref{thm:main}, we
will give a detailed argument that does not quote Schwartz's result.
Of course in many places, our proof parallels some of the intermediate
results in ~\cite{richBook}.

We will not attempt to give an explicit allowable range of parameters
$t$ in Theorem~\ref{thm:main}, although it would certainly be
interesting to do so (and also to try and make this range optimal).

The bulk of the work will be to describe the Ford domain for the
holonomy group of the unipotent uniformization of $M$, and to study in
detail the generic character of the intersection of its sides, along
facets of all dimensions. The genericity that we will prove is
genericity at infinity, namely we will show that each ideal vertex in
the Ford domain lies on precisely three sides that intersect
transversely at that point.  For finite vertices, no genericity is to
be expected, since the group is known to contain elliptic elements of
order 3 and 4 (see~\cite{derauxfalbel}). In fact all the deformations
we consider will preserve the conjugacy classes of these elliptic
elements, and we will show that they do not affect the non-generic
character of the fundamental domains at these points:
\begin{prop}
  The image of $\rho_t$ is a triangle group. More specifically, for
  all $t$, we have
    $$
    \rho_t(g_2)^4=\rho_t(g_1g_2)^3=\rho_t(g_2g_1g_2)^3=id.
    $$
\end{prop}

\begin{flushleft}{\bf Acknowledgements:} 
  This work was partly supported by the ANR, through the grant SGT
  (``Structures G\'eom\'etriques Triangul\'ees''). I also benefited
  from generous support from the GEAR network (NSF grants DMS 1107452,
  1107263, 1107367), via funding of a long term visit at ICERM. I am
  pleased to thank ICERM for its hospitality, and the participants of
  the semester program entitled ``Low-dimensional Topology, Geometry,
  and Dynamics'' for inspiring interactions.  Finally, I would like to
  thank the referee for several suggestions that helped improve the
  readability of the paper.
\end{flushleft}

\section{The real hyperbolic Ford domain}

Throughout this section, we denote by $M$ the figure eight knot
complement. We review the description of a cusp neighborhood for $M$.
This is probably familiar to most readers, but the details will be
used in the identification of the manifold at infinity of our complex
hyperbolic groups. Moreover, quite remarkably, the local combinatorics
of the real hyperbolic Ford domain turn out to be exactly the same as
the local combinatorics of our fundamental domain for the action of
the group on the domain of discontinuity.

Recall that the fundamental group $\pi_1(M)$ has a presentation of the
form
$$
\presGroup{g_1,g_2,g_3}{g_2=[g_3,g_1^{-1}], g_1g_2=g_2g_3},
$$ 
with peripheral subgroup generated by $g_3^{-1}$ and
$g_1(g_1g_2)^{-1}g_3g_2g_3^{-1}$.

From this, one can find all type-preserving representations of
$\pi_1(M)$ up to conjugation, as in~\cite{rileyquadratic}. Indeed, the
generators $g_1$ and $g_3$ should be parabolic elements in $SL_2(\C)$,
which we denote by $G_1$ and $G_3$. We may assume $G_1$ (resp. $G_3$)
fixes $0$ (resp. $\infty$), and since all parabolic elements are
conjugate, we may also assume
$$
G_1=\begin{pmatrix}
1 & 0\\
-\omega & 1
\end{pmatrix}; 
G_3=\begin{pmatrix}
1 & 1\\
0 & 1
\end{pmatrix}
$$ 
for some $\omega\in\C$. The relation
$G_1[G_3,G_1^{-1}]=[G_3,G_1^{-1}]G_3$ in $PSL_2(\C)$, is easily seen
to imply $\omega^2+\omega+1$, so we may take
$$
\omega=\frac{-1+i\sqrt{3}}{2}.
$$

The stabilizer of $\infty$ in $PSL_2(\Z[\omega])$ is clearly given by
translations by Eisenstein integers, but the stabilizer in the group
generated by $G_1$ and $G_3$ is slightly smaller, it can be checked to
be generated by translations by $1$ and $2i\sqrt{3}$
(see~\cite{rileyquadratic} for more details).

\begin{figure}[htbp]
  \epsfig{figure=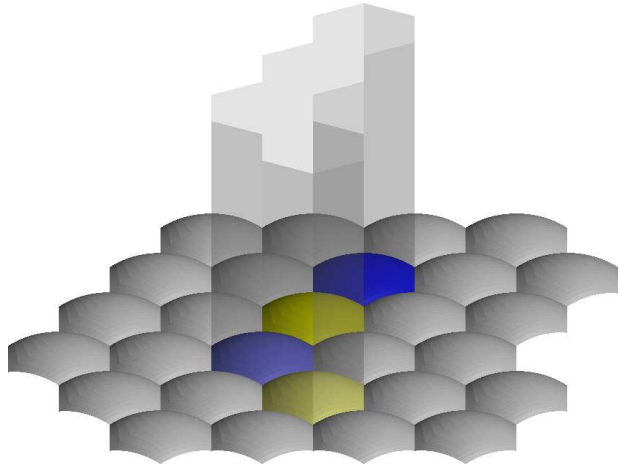,width=0.5\textwidth}
  \caption{A fundamental domain for the action of $\Gamma$ is an
    infinite chimney over the union of four hexagons, each hexagons
    living in a unit hemisphere around the appropriate Eisenstein
    integer.}\label{fig:m004_hyp}
\end{figure}

\begin{figure}[htbp]
\epsfig{figure=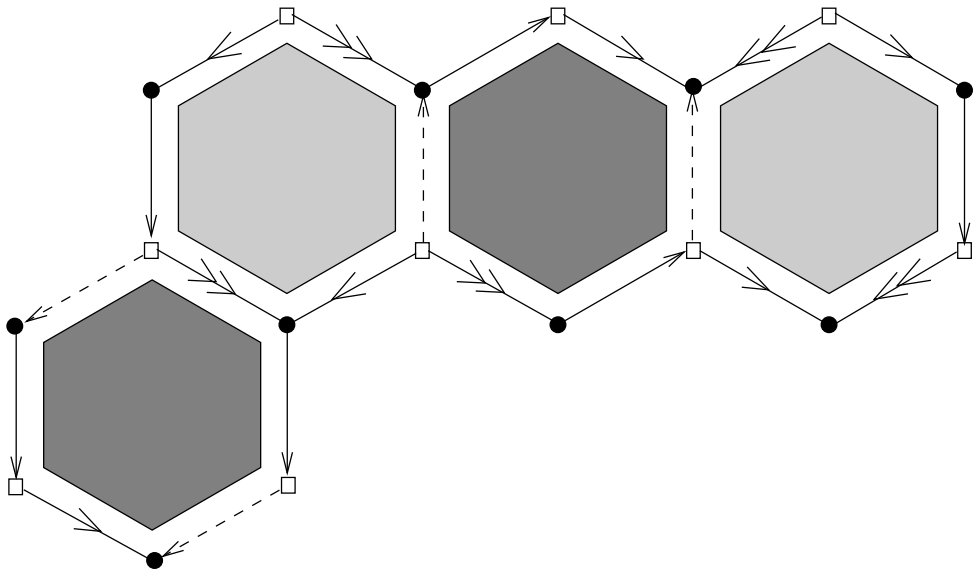, width=0.8\textwidth}
\caption{Bottom of the prism (spine of the figure eight knot
  complement).}\label{fig:torusLink3}
\end{figure}

Recall that the Ford isometric sphere of an element
$$
\begin{pmatrix}
a&b\\
c&d
\end{pmatrix}
$$ 
is bounded by the circle $|cz+d|=1$. The Ford domain turns out to be
the intersection of the exteriors of all spheres of radius 1 centered
at Eisenstein integers.  A schematic picture is shown in
Figure~\ref{fig:m004_hyp}, where the sides corresponding to $G_1^{\pm
  1}$ (resp. $G_2^{\pm 1}=[G_3,G_1^{-1}]^{\pm 1}$) are shaded in the
same color, so the corresponding $2$-faces get identified by the
corresponding isometries.
The complete description of identifications on bottom face of the
prism is given in Figure~\ref{fig:torusLink3}, and there are also
identifications on the vertical sides of the prism, which are simply
given by translations whenever these sides are parallel.
Note that these identifications are described
in~\cite{rileyquadratic}; using current computer technology, they can
also be found using the pictures produced by SnapPy.

\section{Basic complex hyperbolic geometry} \label{sec:chg}

In this section we review some basic material about the complex
hyperbolic plane. The reader can find more details in~\cite{goldman}.

Recall that $\C^{2,1}$ denotes $\C^3$, equipped with a Hermitian form
of signature (2,1). The standard such form is given by $\langle
V,W\rangle = V_1\overline{W}_3 + V_2 \overline{W}_2 + V_3\overline{W}_1 = W^* J V$,
where
$$
J = \left(\begin{matrix}
  0 & 0 & 1\\
  0 & 1 & 0\\
  1 & 0 & 0
\end{matrix}\right).
$$ 
We denote by $U(2,1)$ the subgroup of $GL(3,\C)$ that preserves that
Hermitian form, and by $PU(2,1)$ the same group modulo scalar
matrices. It is sometimes convenient to work with $SU(2,1)$, which is
a 3-fold cover of $PU(2,1)$.

The complex hyperbolic plane $H^2_\C$ is the set of negative complex
lines in $\C^{2,1}$, equipped with a K\"ahler metric that is invariant
under the action of $PU(2,1)$. Such a metric is unique up to scaling,
and it turns out to have constant holomorphic sectional curvature
(which one can choose to be -1).

It is well known that the maximal totally geodesic submanifolds of
$H^2_\C$ are copies of $H^1_\C$ (with curvature -1) and copies of
$H^2_\R$ (with curvature $-1/4$).
 
\subsection{Bisectors} \label{sec:bisectors}

The corresponding distance function is given by
$$
  \cosh^2\frac{1}{2}d(z,w)=\frac{|\langle Z,W\rangle|^2}{\langle Z,Z\rangle \langle W,W\rangle},
$$ 
where $Z$ (resp. $W$) denotes a representative of $z$
(resp. $w$). Given two points $p\neq q\in H^2_\C$, the locus $\B(p,q)$
of points that are equidistant of $p$ and $q$ is called a
bisector. Beware that isometries switching $p$ and $q$ do not fix the
corresponding bisector pointwise, and in fact bisectors are not
totally geodesic. The copies of $\ch 1$ (resp. $\rh 2$) in $\B(p,q)$
are called its complex (resp. real) slices. All real slices intersect
along the same real geodesic, called the {\bf real spine} of the
bisector (see~\cite{goldman}).

Every bisector in $H^2_\C$ is diffeomorphic the unit ball in $\R^3$,
in such a way that the vertical axis is the real spine, complex slices
are horizontal disks, and real slices are disks in vertical planes
containing the vertical axis. One
way to do this explicitly for the bisector $\B(p,q)$ is to scale $q$
by a complex number of modulus one so that $\langle p,q\rangle$ is
real and negative. Then an orthogonal basis for $\C^{2,1}$ is given by
$v_0=p+q$, $v_1=p-q$, $v_2=v_0\boxtimes v_1$ ($\boxtimes$ denotes the
Hermitian cross product, see p.43 of~\cite{goldman}). Of course this
basis can be made Lorentz orthonormal by scaling its vectors so that
$\langle v_0,v_0\rangle =-1$, $\langle v_1,v_1\rangle =1$ and $\langle
v_2,v_2\rangle =1$. The bisector then can be parametrized by
$(z,t)\in\C\times\R$ by taking vectors of the form
$$
v_0+i\,t\,v_1+z\,v_2.
$$

Given a set $S\subset H^2_\C$, we write $\B(S)$ for the locus
equidistant of all point in $S$, which can be thought of as an
intersection of bisectors.

The intersection of two bisectors is usually not totally geodesic, but
it can be in some rare instances. When $p$, $q$ and $r$ are not in a
common complex line (i.e. when lifts of these vectors are linearly
independent), the locus $\B(p,q,r)$ of points equidistant of $p$, $q$
and $r$ is a smooth non totally geodesic disk, and is often called a
Giraud disk, see~\cite{giraud}. The following property is crucial when
studying fundamental domains (see~\cite{giraud},~\cite{goldman}).
\begin{thm}\label{thm:giraud}
  If $p$, $q$ and $r$ are not in a common complex line, then
  $\B(p,q,r)$ is contained in precisely three bisectors, namely
  $\B(p,q)$, $\B(q,r)$ and $\B(q,r)$.
\end{thm}
Note that checking whether an isometry maps a Giraud disk to another
is equivalent to checking that the corresponding triple of points are
mapped to each other.

In order to study Giraud disks, we will use spinal coordinates. The
complex slices of $\B(p,q)$ are given explicitly by choosing a lift
$\tilde{p}$ (resp. $\tilde{q}$) of $p$ (resp. $q$).

When $p,q\in H^2_\C$, we simply choose lifts such that $\langle
\tilde{p},\tilde{p}\rangle = \langle \tilde{q},\tilde{q}\rangle$. In
this paper, we will mainly use these parametrization when
$p,q\in\partial_\infty H^2_\C$. In that case, the condition $\langle
\tilde{p},\tilde{p}\rangle = \langle \tilde{q},\tilde{q}\rangle$ is
vacuous, since all lifts are null vectors;  we then choose some fixed
lift $\tilde{p}$ for the center of the Ford domain, and we take
$\tilde{q}=G\tilde{p}$ for some $G\in U(2,1)$. If a different matrix $G'=SG$,
with $S$ a scalar matrix, note that the diagonal element of $S$ is a
unit complex number, so $\tilde{q}$ is well defined up to a unit
complex number.

The complex slices of $\B(p,q)$ are obtained as (the set of negative
lines in) $(\bar z \tilde{p}-\tilde{q})^\perp$ for some arc of values
of $z\in S^1$, which is determined by requiring that $\langle \bar z
\tilde{p}-\tilde{q}, \bar z \tilde{p}-\tilde{q}\rangle > 0$. 

Since a point of the bisector is on precisely one complex slice, we
can parametrize $\B(p,q,r)$ by $(z_1,z_2)\in S^1\times S^1$ via
\begin{equation}\label{eq:spinalcoords} 
V(z_1,z_2)=(\bar z_1 p-q)\boxtimes(\bar z_2 p - r) = q\boxtimes r +
z_1\ r \boxtimes p + z_2\ p\boxtimes q.
\end{equation}
The Giraud disk corresponds to the $(z_1,z_2)\in S^1\times S^1$ such
that $\langle V(z_1,z_2),V(z_1,z_2)\rangle <0$ (it follows from the
fact that the bisectors are covertical that this region is a
topological disk, but this is not obvious, see chapters~8 and~9
in~\cite{goldman}).

The boundary at infinity $\partial_\infty \B(p,q,r)$ is a circle, given
in spinal coordinates by the equation
\begin{equation}\label{eq:giraudcircle}
\langle V(z_1,z_2),V(z_1,z_2)\rangle = 0
\end{equation}
Note that the choice of two lifts of $q$ and $r$ affects the spinal
coordinates by rotation on each of the $S^1$ factors.

A defining equation for the trace of another bisector $\B(a,b)$ on the
Giraud disk $\B(p,q,r)$ can be written in the form
\begin{equation}\label{eq:thirdbisector}
|\langle V(z_1,z_2),a\rangle| = |\langle V(z_1,z_2),b\rangle|, 
\end{equation}
provided $a,b$ are suitably chosen lifts. The expressions $\langle
V(z_1,z_2),a\rangle$ and $\langle V(z_1,z_2),b\rangle$ are affine in
$z_1$, $z_2$. 

These triple bisector intersections can be parametrized fairly
explicitly, because one can solve the equation
$
|\langle V(z_1,z_2),a\rangle|^2 = |\langle V(z_1,z_2),b\rangle|^2
$
for one of the variables $z_1$ or $z_2$, simply by solving a quadratic
equation. A detailed explanation of how this works can be found in
section~2.3 of~\cite{derauxfalbel}, we will also review this in
section~\ref{sec:rur}.

Note that our parameters also give a parametrization of the
intersection in $P^2_\C$ of the extors extending the bisectors, see
chapter~8 of~\cite{goldman}. The Giraud disk is a disk in the
intersection of the extors, which is a torus.

\subsection{Siegel domain and the Heisenberg group} \label{sec:siegel}

The complex analogue of the upper half space model for $H^n_\R$ is the
Siegel domain, which is obtained by sending the line spanned by
$(1,0,0)$ to infinity. We denote the corresponding point of
$\partial_\infty H^2_\C$ by $p_\infty$.

More precisely, we take affine coordinates $z_1=Z_1/Z_3$,
$z_2=Z_2/Z_3$, and a negative complex line has a unique representative
of the form $z=(z_1,z_2,1)$ with
$$
z^*Jz=2\re(z_1)+|z_2|^2<0
$$

Since we are interested in geometric structures modeled on
$\partial_\infty H^2_\C$, we will use mainly the boundary of the
Siegel domain, which is given by points $z=(z_1,z_2,1)$ with
$2\re(z_1)+|z_2|^2=0$. It is best understood in terms of Heisenberg
geometry, as we now briefly recall.

A large part of the stabilizer of the point at infinity is given by
unipotent upper triangular matrices. One easily checks that such a
matrix preserve the Hermitian form $J$ if and only if it can be
written as
$$
\left(
\begin{matrix}
  1 &  -\bar a\sqrt{2} & -|a|^2+is\\
  0 &   1      &   a\sqrt{2}\\
  0 &   0      &   1
\end{matrix}
\right)
$$ 
for some $(a,s)\in\C\times\R$. Since these upper triangular
matrices form a group, we get a group law on $\C\times \R$, given by
\begin{equation}\label{eq:heis}
(a,s)*(a',s')=(a+a',s+s'+2\Im(a\,\bar a')).
\end{equation}
This is the so-called Heisenberg group law.

The action of the unipotent stabilizer of $p_\infty$ is simply
transitive on $\partial_\infty H^2_\C - \{p_\infty\}$, so we will often
identify the latter with $\C\times \R$.

The boundary at infinity of totally geodesic subspaces can be seen in
somewhat simple terms in $\C\times \R$. The boundary of a copy of
$H^1_\C$ (which is the intersection of an affine line in $\C^2$ with
the Siegel half space) is called a $\C$-circle. These are ellipses
that project to circles in $\C$ (or possibly vertical lines, if they
go through $p_\infty$).

The boundary of copies of $H^2_\R$ (which are images under arbitrary
isometries of the set of real points in the Siegel half space)
intersect the boundary at infinity in a so-called $\R$-circle. In the
Heisenberg group, these are curves that project to lemniscates in $\C$
(or possibly straight lines when they go through $p_\infty$). For more
on this, see chapter~4 of~\cite{goldman}, for instance.

The full stabilizer of $p_\infty$ is generated by the above unipotent
group, together with the isometries of the form 
$$
\left(\begin{matrix}
  1 & 0 & 0\\
  0 & e^{i\theta} & 0\\
  0 & 0 & 1
\end{matrix}\right),\quad 
\left(\begin{matrix}
  \lambda & 0 & 0\\
  0 & 1 & 0\\
  0 & 0 & 1/\lambda
\end{matrix}\right),
$$ 
where $\theta,\lambda\in \R$, $\lambda\neq 0$. The first one acts
on Heisenberg as a rotation with vertical axis:
$$
(a,s)\mapsto (e^{i\theta}a,s),
$$
whereas the second one acts as
$$
(a,s)\mapsto (la,l^2s).
$$

There is a natural invariant metric on the Heisenberg group, called
the Cygan metric, given by $d(g,g')=||g^{-1}g'||$, and the norm of an
element of the Heisenberg group is given by
\begin{equation}\label{eq:cygan}
||(z,t)|| = ||z|^2+it|^{1/2}
\end{equation}
The Cygan sphere with center $(z_0,t_0)$ and radius $r$ has equation
\begin{equation}\label{eq:cygansphere}
  \left||z-z_0|^2+i(t-t_0+2\Im(z\bar z_0))\right|=r^2.
\end{equation}

\subsection{Ford domains and the Poincar\'e polyhedron theorem} \label{sec:defford}

Let $\Gamma$ be a subgroup of $PU(2,1)$, let $q\in\partial_\infty
H^2_\C$ and let $Q$ denote a lift of $q$ in $\C^{2,1}$. 
\begin{dfn}
  The Ford domain for $\Gamma$ centered at $q$ is the set
  $F_{\Gamma,q}$ of points $z\in H^2_\C$ such that
  $$
  |\langle Z,Q \rangle| \leq |\langle Z, G(Q)\rangle|
  $$
  where $G$ is a matrix representative of some element $g\in \Gamma$.
\end{dfn}
The inequality is actually independent of the lift $G\in U(2,1)$
chosen for $g\in PU(2,1)$. For a given $g\in \Gamma$ and lift $G\in
U(2,1)$, we denote by $\B_g$ the bisector given in homogeneous
coordinates by
\begin{equation}\label{eq:Bg}
  |\langle Z,Q \rangle| = |\langle Z, G(Q)\rangle|.
\end{equation}
For concreteness, we mention that the boundary at infinity of $\B_g$
can be described as a Cygan sphere in the Heisenberg group (see
section~\ref{sec:siegel}). The Cygan sphere corresponding to an
element $G$ has radius $\sqrt{2/|g_{31}|}$ (note that $G$ fixes
$p_\infty$ if and only if $g_{31}=0$) and center $(\bar g_{32}/\bar
g_{31},2\Im(\bar g_{33}/\bar g_{31})$ (see equation~\eqref{eq:cygansphere}).

We denote by $b_g=\B_g\cap F$, i.e. the side of $F$ that lies on the
bisector $\B_g$, and we refer to it as {\bf the side corresponding to
  the group element $g$}. For a general $g\in\Gamma$, $b_g$ may have
dimension smaller than 3 (in fact it is often empty). The bisectors of
the form $\B_g$ such that $b_g$ have dimension three will be called
{\bf bounding bisectors}.

The basic fact is that if $q$ has trivial stabilizer in $\Gamma$, then
$F=F_{\Gamma,q}$ is a fundamental domain for its action. However, it
is customary to take $q$ to have a nontrivial stabilizer $H\subset
\Gamma$, in which case $F$ is only a fundamental domain modulo the
action of $H$. In other words, in that case, $F$ is a fundamental
domain for the decomposition of $\Gamma$ into cosets of $H$.

It is ususally very hard to determine $F$ explicitly; in order to
prove that a given polyhedron is equal to $F$, the main tool is the
Poincar\'e polyhedron theorem. The basic idea is that the sides of $F$
should be paired by isometries, and the images of $F$ under these
so-called side-pairing maps should give a local tiling of $H^2_\C$. If
they do (and if the quotient of $F$ by the identifications given by
the side-pairing maps is complete), then the Poincar\'e polyhedron
implies that the images of $F$ actually give a global tiling.

Once a fundamental domain is obtained, one gets an explicit
presentation of $\Gamma$ in terms of the generators given by the
side-pairing maps together with a generating set for the stabilizer
$H$, the relations corresponding to so-called ridge cycles (which
correspond to the local tiling near each codimension two face).

For more details on this theorem, see~\cite{derauxfalbel},~\cite{dpp2}
and~\cite{JRP-book}.

\section{A boundary parabolic family of representations}

In this section, we parametrize a neighborhood of the unipotent
solution in the character variety $\chi(\pi_1(M),PU(2,1))$. We will
use the presentation
$$
\langle\ g_1,g_2,g_3\ |\ g_1g_2=g_2g_3, g_2=[g_3,g_1^{-1}]\ \rangle.
$$ 
In order to describe representations, we seek to parametrize triples
$G_1,G_2,G_3$ of matrices in $SU(2,1)$ that satisfy the same relations
as $g_1$, $g_2$, $g_3$ (possibly up to multiplication by a scalar
matrix, since we are really after representations in $PU(2,1)$).

If the fixed points of $G_1$ and $G_3$ are distinct, we may assume
\begin{equation}\label{eq:genericgroup}
G_1=\left( \begin{matrix} \lambda & a & b\\
                          0 & \overline{\lambda}^2 & c\\
                          0 & 0 & \lambda
\end{matrix}\right),\quad 
G_3=\left( \begin{matrix} \lambda & 0 & 0\\
                          f & \overline{\lambda}^2 & 0\\
                          e & d & \lambda
\end{matrix}\right),
\end{equation}
were $|\lambda|=1$. 

Note that the representation considered in~\cite{derauxfalbel} is obtained by taking
$$
\lambda=1, a=d=1, c=f=-1, b=\overline{e}=-(1+i\sqrt{7})/2
$$
in equation~\eqref{eq:genericgroup}.

The fact that $G_1$ and $G_3$ are isometries of the form $J$ implies
\begin{equation}\label{eq:isom}
\left\{
\begin{array}{l}
  c=-\overline{a}\,\overline{\lambda},\quad f=-\overline{d}\,\overline{\lambda}\\
  |d|^2+ \overline{e} \lambda + e\overline{\lambda}=0\\
  |a|^2+ \overline{b} \lambda + b\overline{\lambda}=0
\end{array}
\right.
\end{equation}
We then compute the commutator $G_2=[G_3,G_1^{-1}]$, and consider the
system of equations given by $R=0$, where
\begin{equation}\label{eq:1223}
  R=G_1G_2-G_2G_3.
\end{equation}
Note that this already restricts the character variety, since we only
consider representations into $U(2,1)$ rather than $PU(2,1)$, but this
is fine if we are after a neighborhood of the boundary unipotent solution,
where the relation~\eqref{eq:1223} holds in $U(2,1)$.

Requiring that $G_1$ and $G_3$ preserve the standard antidiagonal
form, we must have

The $(1,1)$-entry of $R$ is given by
\begin{equation}
(|a|^2 e - |d|^2 b) (1+\overline{a}d-\lambda^3-\overline{\lambda}^3).
\end{equation}
The first factor does not vanish for the boundary unipotent solution,
so in its component we must have
\begin{equation}\label{eq:1-1}
1+\overline{a}d=\lambda^3+\overline{\lambda}^3.
\end{equation}

Note that by conjugation by a diagonal matrix with diagonal entries
$k_1,k_2,k_3$, we can assume that $a\in\R$ (and we can also impose
that $|b|$ is given by any positive real
number). Equation~\eqref{eq:1-1} then implies that $d$ is real as
well, so from this point on we assume
$$
a,d\in\R.
$$

The $(2,2)$-entry of $R$ can then be written as
$$
-(|a|^2 e - |d|^2 b)
    (a^2e\overline{\lambda}^4 + a^2 d^2 \overline{\lambda}^3 - ad + be\overline{\lambda}^5 - 1 +bd^2\overline{\lambda}^4),
$$
so we get the equation
\begin{equation}\label{eq:2-2}
  a^2e\overline{\lambda}^4 + a^2 d^2 \overline{\lambda}^3 - ad + be\overline{\lambda}^5 - 1 +bd^2\overline{\lambda}^4.
\end{equation}
Using the relations~\eqref{eq:isom} and~\eqref{eq:1-1},~\eqref{eq:2-2}, can be rewritten as
\begin{equation}\label{eq:last}
  be\lambda = \lambda^3+\overline{\lambda}^3.
\end{equation}

As mentioned above, by conjugation by a diagonal matrix, we can adjust
$|b|$, for instance so that
$$
|b|^2=\lambda^3+\overline{\lambda}^3,
$$
and in that case~\eqref{eq:last} implies
$$
|e|^2=|b|^2.
$$

We will now show that, given $\lambda$, the following system has
precisely two solutions:
\begin{equation}\label{eq:system}
\left\{\begin{array}{l}
   a^2+\overline{b}\lambda+b\overline{\lambda}=0\\
   d^2+\overline{e}\lambda+e\overline{\lambda}=0\\
   1+ad=\lambda^3+\overline{\lambda}^3\\
   eb\lambda=\lambda^3+\overline{\lambda}^3\\
   |b|^2=\lambda^3+\overline{\lambda}^3.
\end{array}\right.
\end{equation}
In order to do that, note that the first four imply
$$
b\overline{e}+\overline{b}e=1-2(\lambda^3+\overline{\lambda}^3),
$$
and the last two imply 
$$
e=\overline{b}\overline{\lambda}.
$$
Putting these two together, we get
\begin{equation}\label{eq:line-circle}
  \Re(b^2\lambda)=\frac{1}{2}-2\kappa,
\end{equation}
where we have written
\begin{equation}
  \kappa=(\lambda^3+\overline{\lambda}^3)/2.
\end{equation}

The equation $\Re(z)=\frac{1}{2}-2\kappa$ has a solution with $|z|=2\kappa$ if and only if
$$
2\kappa\geq \frac{1}{2}-2\kappa,
$$ 
and in that case one gets a simple formula for the solutions
(intersect a vertical line with the circle of radius $|2\kappa|$ centered
at the origin).

We get that~\eqref{eq:line-circle} has solutions if and only if $\kappa\geq
\frac{1}{8}$, and the solutions are given by
\begin{equation}\label{eq:b2l}
  b^2\lambda = \frac{1}{2}-2\kappa \pm i \sqrt{\frac{1}{2}(4\kappa-\frac{1}{2})}.
\end{equation}
This determines $b$ up to its sign, opposite values clearly giving
conjugate groups (they differ by conjugation by a diagonal
matrix). The two values also yield isomorphic groups, obtained from
each other by complex conjugation.

We will choose the solution to match the notation for the unipotent
solution given in~\cite{derauxfalbel}, which corresponds to
$\lambda=1$, $a=d=1$, $b=-\frac{1+i\sqrt{7}}{2}$ and
$e=-\frac{1-i\sqrt{7}}{2}$.

As a consequence, we take
$$
b=-\frac{1+i\sqrt{8\kappa-1}}{2\sqrt{\lambda}},
$$
where we take the squareroot to vary continuously near $\lambda=1$.

The system~\eqref{eq:system} then gives values for the other
parameters, namely
$$
e=2\kappa/b\lambda=-\frac{1-i\sqrt{8\kappa-1}}{2\sqrt{\lambda}},
$$ 
and one easily writes an explicit formula for $a$ and $d$ (once
again, these are determined only up to sign, but changing $a$ to $-a$
can be effected by conjugation by a diagonal matrix). The formula is
as follows,
$$
a=\sqrt{(4\mu^2-3)\mu + \sqrt{8\kappa-1} (4\mu^2-1)\nu}, \quad d=\sqrt{(4\mu^2-3)\mu - \sqrt{8\kappa-1} (4\mu^2-1)\nu},
$$ 
where we have written $\sqrt{\lambda}=\mu+i\nu$ with $\mu$, $\nu$
real. In terms of this new parameter, the condition $\kappa>1/8$
translates into
$$
\mu>\cos(\frac{1}{3}\arctan\frac{\sqrt{7}}{3}) = 0.9711209254\dots
$$
In fact, in order to get $a$ and $d$ to be real, we also need
$$
(4\mu^2-3)\mu - \sqrt{8\kappa-1} (4\mu^2-1)\nu\geq 0,
$$ 
which translates into $\mu\geq \cos(\pi/18)$. The value
$\mu=\cos(\pi/18)$ corresponds to a situation where $d=0$.

\subsection{Triangle group relations}

The following matrices can be computed explicitly:
$$
G_2 = \left(\begin{matrix}
  1+\l^3       &  a\lc-\overline{b}d    & (e+b)\lc\\
  ab-d\lc^2    &    -\l^3               & 0\\
  (e+b)\lc     &      0                 & 0
\end{matrix}
\right)
$$
$$
G_1G_2 = \left(\begin{matrix}
  \l                   &  a(1-\l^3)-ed\l^2          & (e+b)\\
  -\lc^2(ae+d\lc^2)    &    -\l                     & 0\\
  (e+b)                &      0                     & 0
\end{matrix}
\right)
$$
$$
G_1^2G_2 = \left(\begin{matrix}
  \lc          &  -\l^3(a\l+ed)         & (e+b)\l\\
  \l^2(ab+d\l) &    -\lc                & 0\\
  (e+b)\l      &      0                 & 0
\end{matrix}
\right)
$$

In particular, 
$$
{\rm tr}(G_2)=1,\quad {\rm tr}(G_1G_2)=0,\quad {\rm tr}(G_2G_1G_2)=0,
$$ 
or in other words,
$$
G_2^4=id,\quad (G_1G_2)^3=id,\quad (G_1^2G_2)^3=id.
$$
The last two relations imply that
$$
(G_2G_1G_2)^3=id.
$$

\begin{prop}\label{prop:twistrelations}
  Throughout the twist parabolic deformation, we have $G_1G_2=G_2G_3$,
  $G_2=[G_3,G_1^{-1}]$, $G_2^4=id$, $(G_1G_2)^3=id$,
  $(G_2G_1G_2)^3=id$.
\end{prop}

\subsection{Fixed points of elliptic elements}\label{sec:fixedpoints}

Note also that for each of the three matrices $G_2$, $G_1G_2$ and
$G_1^2G_2$, the negative eigenvector is the one with eigenvalue $1$
(indeed, this is true for the unipotent solution, so it holds
throughout the corresponding component of the character variety).

For future reference, we give explicit formulas for these fixed
points:
$$
p_{2} = \left(1+\l^3, ab-d\lc^2, (\lc+\l^2)(e+b)\right),
$$
$$
p_{12} = \left(1+\l, -\lc^2(ae+d\lc^2), (1+\lc)(e+b)\right),
$$
$$
p_{112} = \left(1+\lc, \l^2(ab+d\lc), (\lc+\lc^2)(e+b)\right).
$$

\begin{lem}\label{lem:p2}
  Throughout the deformation, $p_2$ is on six bounding bisectors,
  corresponding to the group following elements
  $$ 
  {2},\ {\bar 2},\ {3},\ {12},\ {\bar 1\bar 2},\ {\bar 1 3}.
  $$
\end{lem}
\begin{pf}
  The statement about $G_2^{\pm 1}$ is obvious since $p_2$ is fixed by
  $G_2$. The other four statements all follow from 
  \begin{equation}\label{eq:key}
    d(p_2,p_0)=d(p_2,(G_2G_1)^{-1}p_0).
  \end{equation}
  Indeed,
  $$
  d(p_2,(G_2G_1)^{-1}p_0)=d(p_2,G_2^{-1}G_1^{-1}G_2^{-1} p_0)=d(p_2,G_1G_2 p_0),
  $$
  where we have use $G_1p_0=p_0$ and $(G_1G_2)^3=id$. Similarly,
  using $G_1G_2=G_2G_3$, we get
  $$
  d(p_2,G_1G_2p_0)=d(p_2,G_2^{-1}G_1G_2 p_0)=d(p_2,G_3 p_0).
  $$
  Finally, using $G_2=[G_3,G_1^{-1}]$ we get
  $$
  d(p_2,G_3p_0)=d(p_2,G_2^{-1}G_3 p_0)=d(p_2,G_1^{-1}G_3 p_0).
  $$

  In order to prove~\eqref{eq:key}, we compute
  $$
  G_1^{-1}G_2^{-1} p_0 = (\overline{b}+\overline{e})\l \left( \overline{b}, a , \lc \right),
  $$
  and we observe $|(\overline{b}+\overline{e})\l|=1$, so we need only check
  $$
  |\langle p_2, p_0 \rangle|=|\langle p_2, X\rangle|,
  $$
  where $X=\left( \overline{b}, a , \lc \right)$.
  Now 
  $$
  |\langle p_2, p_0 \rangle|^2 = |(\l+\lc^2)(\overline{e}+\overline{b})|^2=|1+\lambda^3|^2=2+\l^3+\lc^3, 
  $$
  and
  $$
  \langle p_2, X\rangle = \lc (2-\l^3-\lc^3-b^2\l), \qquad  |\langle p_2, X\rangle|^2=2+ \l^3+\lc^3.
  $$
\end{pf}

\begin{lem}\label{lem:p1b21}
  Through the deformation, $p_{\bar121}=G_1^{-1}p_2$ stays on six bounding
  bisectors, corresponding to the following group elements:
  $$ 
  2,\ \bar 1 2,\ \bar1\bar2,\ \bar1 3,\ \bar1\bar1\bar 2,\ \bar1\bar13.
  $$
\end{lem}

\begin{pf}
  The statement follows from Lemma~\ref{lem:p2} by conjugation by
  $G_1^{-1}$ (which by definition fixes $p_0$).
\end{pf}

\section{Combinatorics of the Ford domain in the unipotent case}

In this section, we denote by $\Gamma$ the image of $\rho_0$. It is
generated by the matrices
$$
G_1=\begin{pmatrix}
           1&  1 &   -\frac{1}{2}- \frac{\sqrt{7}}{2}i\\
 0 & 1 & -1  \\
0 & 0 & 1\\
           \end{pmatrix},
 \ \
            G_3=\begin{pmatrix}
           1&  0 &   0\\
 -1 & 1 & 0  \\
-\frac{1}{2}+ \frac{\sqrt{7}}{2}i& 1 & 1\\
           \end{pmatrix}.
$$ 
One then sets
$$
G_2=[G_3,G_1^{-1}].
$$ 
We will often use word notation in the generating set $G_1$, $G_2$,
$G_3$, using bars to denote inverses. For instance, $23\bar13$ denotes
$G_2G_3G_1^{-1}G_3$.

We consider the Ford domain centered at the fixed point of $G_1$,
which is $p_\infty$ in the notation of section~\ref{sec:defford}, and
work in the Siegel half space. We denote by $P=\langle G_1\rangle$,
and by $F$ the corresponding Ford domain. We wish to prove that $F$ is
a fundamental domain for the action of the cosets of $P$ in $\Gamma$.

We denote by $S=\{G_2,G_2^{-1},G_3,G_3^{-1}\}$, and by $S^P$ the set
of all conjugates of elements of $S$ by powers of $G_1$. We consider
the partial Ford domain $D$ defined in homogeneous coordinates $Z$
by the inequalities
  $$
  |\langle Z,Q \rangle| \leq |\langle Z, G(Q)\rangle|
  $$ 
for all $G\in S^P$. Clearly $F\subset D$, but we mean to prove:
\begin{thm}\label{thm:partialford}
  $F=D$.
\end{thm}

The key steps in the proof of Theorem~\ref{thm:partialford} will be the following:
\begin{itemize}
  \item Determine the combinatorics of $D$;
  \item Show that the elements in $S^P$ define side-pairing maps for
    $D$;
  \item Verify the hypotheses of the Poincar\'e polyhedron theorem.
\end{itemize}

\subsection{Statement of the combinatorics} \label{sec:statecombo}

Clearly $D$ is $G_1$-invariant, so it is enough to describe the
combinatorics of the sides corresponding to $g\in S$,
i.e. $g=G_2,G_3,G_2^{-1},G_3^{-1}$. We will call the corresponding four
sides $b_1$, $b_2$, $b_3$ and $b_4$ respectively, and refer to them as
core sides; the corresponding bisectors will be denoted by $\B_1$,
$\B_2$, $\B_3$ and $\B_4$. The spinal spheres at infinity of these
four bisectors will be denotes by $\Ss_1$, $\Ss_2$, $\Ss_3$, $\Ss_4$.

We will sometimes index other sides than the four basic sides just
described, mostly when describing computations that would unreasonable
to perform by hand. We will order them by concatenating sets of four
conjugates of the base group elements $2,-2,3,-2$ by different powers
of $G_1$, powers being arranged by increasing values of the absolute
values of the exponent (positive powers first). The words
corresponding to the first 20 bisectors are given by
\begin{eqnarray*}
&2, \bar 2, 3, \bar 3, 
\quad 12\bar1, 1\bar 2\bar 1, 13\bar1, 1\bar3\bar1,
\quad \bar121, \bar1\bar21, \bar131, \bar1\bar31,\qquad\qquad\qquad\\
&\quad 1^22\bar1^2, 1^2\bar2\bar1^2, 1^23\bar 1^2, 1^2\bar3\bar1^2,
\quad \bar1^221^2, \bar1^2\bar21^2, \bar1^231^2, \bar1^2\bar31^2,\quad \dots
\end{eqnarray*}
For example, $\B_{5}=G_1(\B_1)$ is the bisector corresponding to
$G_1G_2G_1^{-1}$ (or equivalently for $G_1G_2$, since $G_1$ fixes the
center of our Ford domain), $\B_{10}=G_1^{-1}(\B_2)$ is the bisector
for $G_1^{-1}G_2^{-1}G_1$.

\begin{figure}[htbp]
  \subfigure[$G_2$]{
    \epsfig{figure=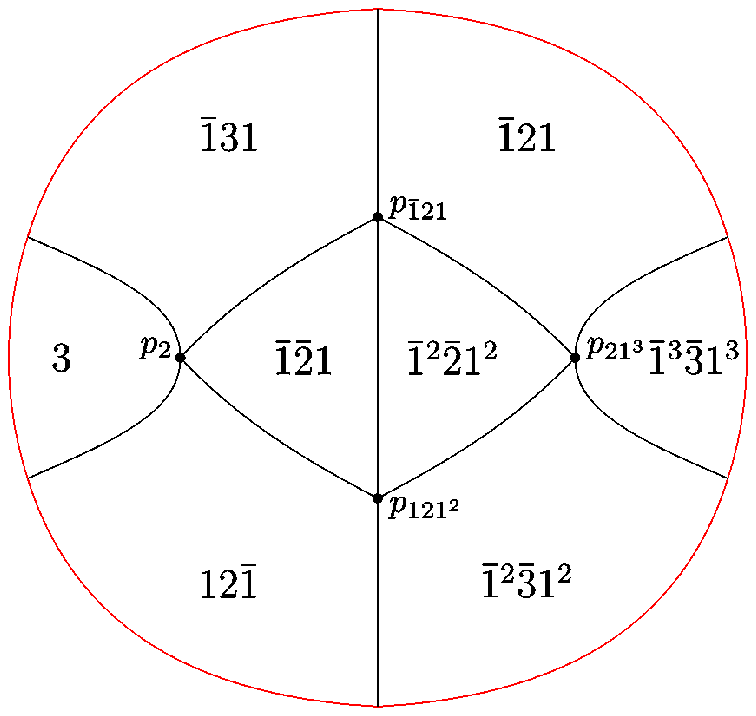, width=0.48\textwidth}
  }\hfill
  \subfigure[$G_2^{-1}$]{\epsfig{figure=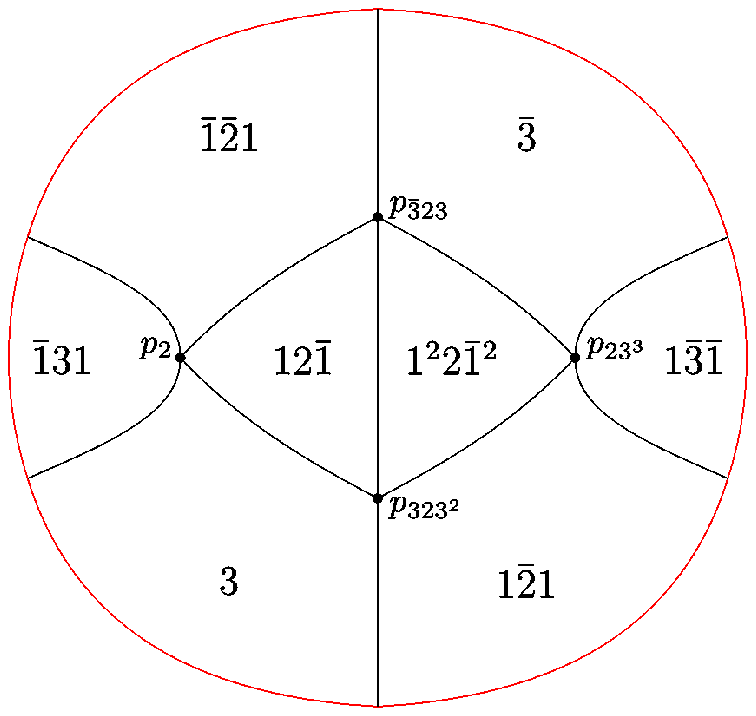, width=0.48\textwidth}}
  \caption{The combinatorics of the face corresponding to $G_2$ and
    $G_2^{-1}$; all 2-faces are labelled, except for the boundary at
    infinity, which is a disk bounded by the most exterior curve
    (shown in red). We also label the finite vertices, namely for
    $w\in \Gamma$, $p_w$ denotes being the isolated fixed point of the
    group element corresponding to the word $w$ ($1=G_1$, $2=G_2$,
    $3=G_3$, $\bar1=G_1^{-1}$, etc).}\label{fig:face2}
\end{figure}

\begin{figure}[htbp]
  \subfigure[$G_3$]{\epsfig{figure=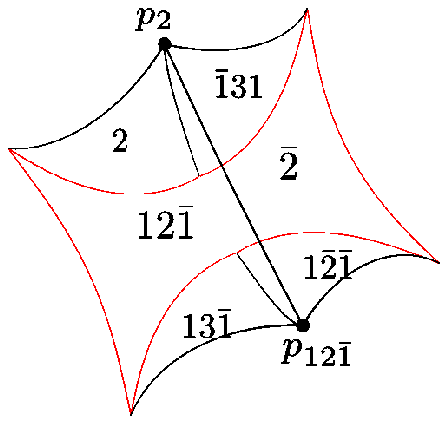, width=0.48\textwidth}}\hfill
  \subfigure[$G_3^{-1}$]{\epsfig{figure=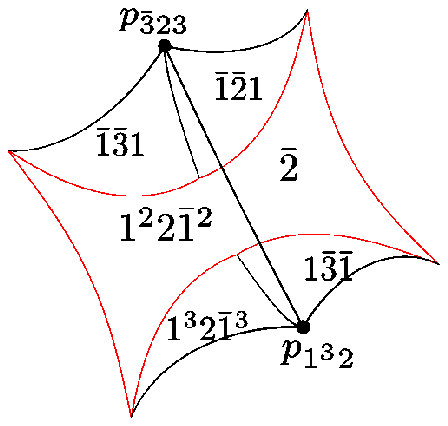, width=0.48\textwidth}}
  \caption{The combinatorics of the face corresponding to $G_3$ and
    $G_3^{-1}$.}\label{fig:face3}
\end{figure}

We describe their combinatorics in the form of pictures, see
Figures~\ref{fig:face2}-\ref{fig:face3}. Each picture is drawn in
projection from a picture where the bisector is identified with the
unit ball in $\R^3$ (see section~\ref{sec:bisectors}). Concretely, we
use spinal coordinates on 2-faces, and parametrize 1-faces by solving
equations of the form~\eqref{eq:thirdbisector} for one of the
variables. 

We also give a list of vertices on the core sides, and also a list of
the bounding bisectors that each vertex lies on,
see Tables~\ref{tab:vertices2} and~\ref{tab:vertices3}.
\begin{table}[htbp]
  \begin{tabular}{c|c|c}
     Word         & bounding bisectors                                                                      & Indices \\ \hline
    $2$          & $2, \bar2,   3,           12\bar1,             \bar1\bar2 1,    \bar131$                & 1,2,3,5,10,11\\
    $\bar1 21$   & $2, \bar121, \bar1\bar21,     \bar131,             \bar1^2\bar21^2, \bar1^231^2$        & 1,9,10,11,18,19\\
    $21^3$       & $2, \bar121, \bar1^2\bar21^2, \bar1^2\bar31^2, \bar1^3\bar21^3, \bar1^3\bar31^3$        & 1,9,18,20,26,28\\
    $121^2$       & $2, 12\bar1, \bar1\bar21,     \bar1\bar31,         \bar1^2\bar21^2, \bar1^2\bar31^2$   & 1,5,10,12,18,20
  \end{tabular}
  \caption{Finite vertices on the face for $G_2$. For each vertex $v$,
    we give a word $w$ for an element that fixes precisely $v$, and
    list the words for the bounding bisectors that contain $v$.}\label{tab:vertices2}
\end{table}
\begin{table}[htbp]
  \begin{tabular}{c|c|c}
    Word         & bounding bisectors                                                                      & Indices\\\hline
    $2$          & $2,     \bar2,   3,           12\bar1,             \bar1\bar2 1,    \bar131$            & 1,2,3,5,10,11\\
    $\bar3 23$   & $\bar2, \bar3,   12\bar1,         \bar1\bar21,         \bar1\bar31,     1^22\bar1^2$    & 2,4,5,10,12,13\\
    $23^3$       & $\bar2, \bar3,   1\bar2\bar1,     1\bar3\bar1,         1^22\bar1^2,     1^32\bar1^3$    & 2,4,6,8,13,21\\
    $323^2$       & $\bar2, 3,       12\bar1,         1\bar2\bar1,         13\bar1,          1^22\bar1^2$  & 2,3,5,6,7,13
  \end{tabular}
  \caption{Finite vertices on the face for $G_2^{-1}$.}\label{tab:vertices3}
\end{table}

\begin{figure}[htbp]
\epsfig{figure=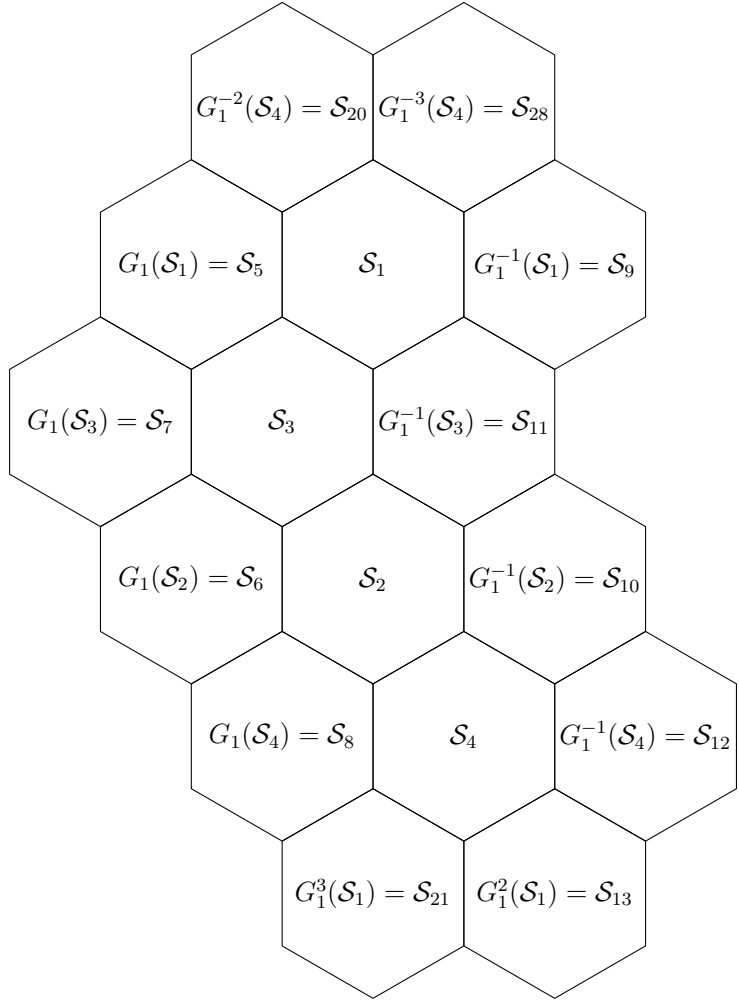, width=0.6\textwidth}
\caption{The combinatorics at infinity of the fundamental domain, near
  the faces for $G_2^{\pm}$ and $G_3^{\pm}$, which are representatives
  of all faces modulo the action of $G_1$.}\label{fig:hexagons}
\end{figure}

\subsection{Effective local finiteness} \label{sec:finite}

The goal of this section is to show that a given face of the Ford
domain intersects only finitely many faces. Since the domain is by
construction $G_1$-invariant, we start by normalizing $G_1$ in a
convenient form.
We will work in the Siegel half space, see section~\ref{sec:siegel}.

A natural set of coordinates is obtained by arranging that $G_2^2$ maps
$p_\infty$ to the origin in the Heisenberg group. There is a unique
Heisenberg translation that achieves this, given by 
$$
Q = \left(
\begin{matrix}
  1 & \frac{3-i\sqrt{7}}{4}  & -\frac{1}{2}\\
  0 &          1             & \frac{-3-i\sqrt{7}}{4}\\
  0 &          0             &         1
\end{matrix}
\right).
$$
One then gets
$$
QG_1Q^{-1} = \left(
\begin{matrix}
  1 & 1 & -\frac{1}{2}\\
  0 & 1 & -1\\
  0 & 0 & 1
\end{matrix}
\right),
$$
and
$$
QG_2^2Q^{-1} = \left(
\begin{matrix}
  0 & 0 & -\frac{1}{2}\\
  0 & -1 & 0\\
  -2 & 0 & 0
\end{matrix}
\right).
$$ 
Of course one could make the last matrix even simpler by composing
with a loxodromic element.

We denote by $A_j=QG_jQ^{-1}$. We then have
$$
A_2(\infty)=(\alpha,0),\quad A_2^2(\infty)=(0,0),\quad A_2^{-1}(\infty)=(-\alpha,0)
$$
where $\alpha=\frac{3+i\sqrt{7}}{4\sqrt{2}}$.
$$ 
A_3(\infty)=(-\frac{1}{2\sqrt{2}},-\frac{\sqrt{7}}{8}),\quad
A_3^{-1}(\infty)=(-\frac{1}{\sqrt{2}},\frac{\sqrt{7}}{2}).
$$

The spinal sphere with center $(0,0)$ and radius $r$ has
equation
\begin{equation}\label{eq:unitsphere}
(x^2+y^2)^2+t^2=r^4,
\end{equation}
so we get a spinal sphere centered at $(a+ib,u)$ by translation:
\begin{equation}\label{eq:std}
\left((x-a)^2+(y-b)^2\right)^2+(t-u-ay-bx)^2=r^4.
\end{equation}
By writing out the equation~\eqref{eq:Bg}, squaring both sides and
identifying with equation~\eqref{eq:std}, one checks that the spheres
$\Ss_1$, $\Ss_2$ have radius $1$, whereas $\Ss_3$, $\Ss_4$ have radius
$2^{-1/4}$. We summarize this information in Table~\ref{tab:corespheres}.
\begin{table}
  \begin{tabular}{ccc}
    Sphere & Center & radius\\
    $\Ss_1$ & $(\frac{3+i\sqrt{7}}{4\sqrt{2}},0)$ & $1$\\
    $\Ss_2$ & $(-\frac{3+i\sqrt{7}}{4\sqrt{2}},0)$ & $1$\\
    $\Ss_3$ & $(-\frac{1}{2\sqrt{2}},-\frac{\sqrt{7}}{8})$ & $2^{-1/4}$\\
    $\Ss_4$ & $(-\frac{1}{\sqrt{2}},\frac{\sqrt{7}}{2})$ & $2^{-1/4}$\\
  \end{tabular}
  \caption{Centers and radii of core spinal spheres.}\label{tab:corespheres}
\end{table}

The action of $A_1$ on the Heisenberg group is given by 
\begin{equation}\label{eq:tsl}
  (z,t)\mapsto (z-1,t+\Im(z)),
\end{equation}
and in particular we get the following:
\begin{prop}
  The element $A_1$ preserves every $\R$-circle of the form
  $(x,0,t_0)$, $x\in\R$. 
\end{prop}
Recall that $\R$-circles are by definition given by the trace at
infinity of totally geodesic copies of $H^2_\R$ in $H^2_\C$. The
corresponding real planes in $H^2_\C$ are preserved by $A_1$, and
their union is the so-called \emph{invariant fan} of $A_1$
(see~\cite{goldmanparkerdirichlet}).

Among all these $\R$-circles, the $x$-axis is somewhat special because
of the following:
\begin{prop}
  The $\R$-plane bounded by the $x$-axis contains the fixed point of
  $G_2$.
\end{prop}
Indeed, the fixed point of $A_2$ is given by
$$
V=(-\frac{1}{2},0,1),
$$
and for $W=(-x^2+it,x\sqrt{2},1)$, 
$$
\langle V,p_\infty\rangle  \langle p_\infty,W \rangle  \langle W, V\rangle = -\frac{1}{2} (1+x^2) + it
$$ 
which is real if and only if $t=0$.

Note that equation~\eqref{eq:tsl} shows that for any two bisectors
$\B_1$ and $\B_2$ not containing $p_\infty$,
$G_1^k\B_1\cap\B_2=\emptyset$ whenever $k$ is large enough. Indeed, it
follows from the detailed study of bisector intersection
in~\cite{goldman} that, if two bisectors intersect, then the
corresponding spinal spheres must intersect.

Moreover, this claim can easily be made effective, i.e. one can get
explicit bounds on how large $k$ needs to be for the above
intersection to be empty.  If $\Ss_j=\partial_\infty\B_j$ is
contained in a strip $\alpha_j\leq x \leq \beta_j$, one can simply
take $k>\beta_2-\alpha_1$, or $k<\alpha_2-\beta_2$.
Note that bounds $\alpha_j,\beta_j$ can be computed fairly easily from
the equations of the relevant spinal spheres (see the
Table~\ref{tab:corespheres} giving the centers and radii).

In particular, we get:
\begin{prop}\label{prop:finite}
  \item $\Ss_1$ intersects $G_1^k\Ss_1$ only if $-2\leq k\leq 2$;\quad $\Ss_1$ intersects $G_1^k\Ss_2$ only if $-4\leq k\leq 1$;
  \item $\Ss_1$ intersects $G_1^k\Ss_3$ only if $-3\leq k\leq 1$;\quad $\Ss_1$ intersects $G_1^k\Ss_4$ only if $-4\leq k\leq 0$.
  \item $\Ss_2$ intersects $G_1^k\Ss_2$ only if $-2\leq k\leq 2$;\quad $\Ss_2$ intersects $G_1^k\Ss_3$ only if $-2\leq k\leq 2$;
  \item $\Ss_2$ intersects $G_1^k\Ss_4$ only if $-2\leq k\leq 2$;\quad $\Ss_3$ intersects $G_1^k\Ss_3$ only if $-2\leq k\leq 2$;
  \item $\Ss_3$ intersects $G_1^k\Ss_4$ only if $-2\leq k\leq 1$;\quad $\Ss_4$ intersects $G_1^k\Ss_4$ only if $-2\leq k\leq 2$;
\end{prop}
This is not an optimal result, since it takes into account only the
variable $x$ and the fact that $G_1$ translates by one unit in the
direction of the $x$-axis. The optimal result is not far from this
though, the point of Proposition~\ref{prop:finite} is to get down to a
finite list of bounding bisectors intersecting a given one (so that we
can use effective computational tools). We will give much more precise
information in the next section.

\subsection{Proof of the combinatorics} \label{sec:provecombo}

The techniques we use in order to justify the combinatorics are very
similar to the ones explained in detail in~\cite{derauxfalbel}
and~\cite{dpp2}. Note that one can think of justifying the
combinatorics as a special case of finding the connected components of
(many) semi-algebraic sets. Indeed, $F$ is clearly semi-algebraic,
defined by inequalities, indexed by $I=\N$:
$$
F=\{ z\in \C^2\ :\forall i\in I,\ f_i(z)<0\},
$$ 
For convenience, we make the convention that $f_0(z)<0$ is the
defining equation for the unit ball, in other words
$$
f_0(z)=\langle \tilde{z}, \tilde{z}\rangle,
$$ 
where $\tilde{z}=(z,1)$. In particular, we consider the boundary at
infinity of complex hyperbolic space as a bounding face.  All other
equations have the form $f_j<0$ where
$$
f_j(z)=|\langle \tilde{z}, \tilde{p}_0\rangle|^2 - |\langle \tilde{z}, \gamma_j \tilde{p}_0\rangle|^2.
$$

The {\bf facets} are of $F$ described by taking some subset $J\subset I$,
and replacing the inequalities indexed by elements of $J$ by the
corresponding equality:
$$
F_J=\{ z\in \C^2\ :\ \forall j\in J,\ f_j(z)=0, \forall i\in I\setminus J,\ f_i(z)<0 \}.
$$ 
The fact that $I$ is infinite will not be a problem because of the
results in section~\ref{sec:finite}, which imply that our polytope is
be locally finite.

More generally, we will consider sets of the form
$$
  F_{J,K}=\{ z\in \C^2\ :\ \forall j\in J,\ f_j(z)=0, \forall i\in K,\ f_i(z)<0 \},
$$ 
where $J$ and $K$ are disjoint. In particular $F_J$ is the same as
$F_{J,I\setminus J}$, and $F_{J,\emptyset}$ is the $|J|$-fold bisector
intersection containing $F_J$.

\subsubsection{Terminology and specification} \label{sec:spec}

We will call {\bf $k$-faces} the facets of our polytopes that have
dimension $k$. Moreover, 3-faces will be simply called {\bf sides},
2-faces will be called {\bf ridges}, 1-faces will be called {\bf
  edges}, and 0-faces will be called {\bf vertices}.

In terms of computations, it will be important to encode vertices.
These can be of two kinds, namely they can be of the form
$F_{A,\emptyset}$, for some $A$ with $|A|=4$, or they can be singular
points of $F_{B,\emptyset}$ with $|B|=3$). In both cases, they can be
obtained by solving a 0-dimensional system (this is the content of
assumption~\ref{assumption1-2}). For each of them, we encode the vertex by
storing a rational univariate representation for the corresponding
solution set, and an isolating interval specifying a root of the
rational parameter (see section~\ref{sec:rur}). 

Note that in the above description, the set $A$ is not unique, since a
vertex may in general lie on more than four bisectors (see the
discussion in section~\ref{sec:fixedpoints}, where we saw examples of
vertices lying on at least six bounding bisectors). Moreover, in
general one cannot take $A$ to be just any $4$-tuple of bisectors that
contain that vertex, since some intersections may not be generic.

We will also need to encode 1-faces. There are two kinds of 1-faces,
namely those that lie in triple bisector intersections (we call these
finite 1-faces), and those that lie in the intersection of the sphere
at infinity $\partial_\infty \ch 2$ with the closure in $\chb 2$ of a
bisector intersection (we call these ideal 1-faces, or 1-faces at
infinity). Computationally, we make no difference between these two
kinds of 1-faces, since both kinds are given in terms of spinal
coordinates for a bisector intersection by an equation that is
quadratic in both variables.

We call an \emph{arc} a subset in $\chb 2$ of a triple bisector
intersection (or a subset of the trace at infinity of a double
bisector intersection) that is
\begin{itemize}
\item homeomorphic to a closed interval,
\item parametrized by one of the spinal coordinates, and such that
\item its endpoints are vertices of the polytope, but its interior
  contains no vertex of the polytope.
\end{itemize}
Note that a 1-face can always be described as a union of finitely many
arcs (but one arc may not suffice, think of a polytope that has a
whole Giraud disk as a facet, so that the boundary of that Giraud disk
is a 1-face homeomorphic to a circle).

We now expand a little on how to parametrize (pieces of) 1-faces by a
single coordinate (we discuss only parametrization by $t_1$, the other
one is entirely similar). Recall from section~\ref{sec:bisectors} that
the relevant defining functions $h(t_1,t_2)$ for triple bisector
intersections (or trace at infinity of double bisector intersections)
have degree at most two in each variable, so we can write them as
$$
a_2(t_1)t_2^2+a_1(t_1)t_2+a_0(t_1),
$$ 
with $a_j$ at most quadratic. With respect to projection onto the
first coordinate axis, the curve usually has two branches, given by
$$
t_2=\frac{-a_1(t_1)\pm \sqrt{\Delta(t_1)}}{2a_2(t_1)},
$$
where
$$
\Delta(t_1)=a_1(t_1)^2-4a_2(t_1)a_0(t_1).
$$ 
Specifically, this occurs above intervals of $t_1$ such that
$a_2(t_1)$ does not vanish. Above such an interval, the ``top branch''
is obtained by taking $+\sqrt{\Delta}$ when $a_2(t_1)>0$, and
$-\sqrt{\Delta}$ when $a_2(t_1)>0$. We call the other branch the
``bottom branch''.

If $a_2$ is identically zero, then the curve is either empty, or
consists of a single vertical line (so branches above the $t_1$ axes
are undefined, and there is a single branch with respect to the
projection onto the $t_2$ axis).

If $a_2$ is not identically zero, it vanishes at one or two points,
and above each of these points, one can determine check whether the
curve contains one, two or infinitely many points (one needs to
determine whether $a_1$, $a_0$ also vanish at these points).

\subsubsection{General procedure} \label{sec:algo}

The pictures in Section~\ref{sec:statecombo} include the statement
that each facet is topologically (in fact piecewise smoothly) a disk
with piecewise smooth boundary (with pieces of the boundary
corresponding to facets of codimension one higher). This is not at all
obvious; one of the difficulties is the fact that the sets $F_{J}$ are
in general not connected, in strong contrast with Dirichlet or Ford
domains in the context of constant curvature geometries (see the
discussion in~\cite{deraux4445}).

For a given $J$, $K$, there is an algorithm to decide whether
$F_{J,K}$ is empty or not, and furthermore one can list its connected
components (and even produce triangulations). One possible approach to
this is the cylindrical algebraic decomposition of semi-algebraic
sets, see~\cite{semialgebraic} for instance.

The main issue when using such algorithms is that the number of
semi-algebraic sets to study is extremely large. If $F$ has $N$ faces,
in principle one has to deal with $N \choose k$ potential facets of
codimension $k$, where $k=1,2,3,4$, which is a fairly large number of
cylindrical decompositions.  Rather, we will bypass the cylindrical
decomposition and use as much geometric information as we can in order
to restrict the number of verifications. Also, rather than using
affine coordinates in $\C^2$, we use natural parametrizations for
bisector intersections, deduced from spinal coordinates (see
section~\ref{sec:bisectors}).

Going back to geometry, the inequality defining complex hyperbolic
space in $\C^2$ (which corresponds to $f_0$) is of course a bit
different from the other inequalities. In particular, when using the
notation $F_{J,K}$, we will always assume one of the index sets $J$ or
$K$ contains $0$.

If $K$ contains $0$, then by definition $F_{J,K}$ is contained in $\ch
2$; we will denote by $\widehat{F}_{J,K}$ its extension to projective
space, namely
$$
\widehat{F}_{J,K} = F_{J,K\setminus \{0\}}.
$$ 
We will also refer to the following set as the trace at infinity of
$F_{J,K}$,
$$
\partial_\infty F_{J,K} = F_{J\cup\{0\},K\setminus \{0\}}.
$$
By $\overline{F}_{J,K}$, we mean the set obtained from the definition of $F_{J,K}$ by replacing $<$ by $\leq$, i.e.
$$ 
\overline{F}_{J,K}=\{ z\in \C^2\ :\ \forall j\in J,\ f_j(z)=0, \forall
i\in K,\ f_i(z)\leq 0 \},
$$ 
which is also
$$
\overline{F}_{J,K} = \bigcup_{L\subset K} F_{J\cup L,K\setminus L}.
$$
Note that in general, this is not the closure of $F_{J,K}$ in $\C^2$.

We focus on an algorithm for determining the combinatorics of ridges,
or in other words facets of the form $F_J$ with $|J|=2$. In most
cases, we will also assume $0\notin J$, i.e. we study finite facets
rather than faces in $\partial_\infty \ch 2$.  The algorithm will
produce a description of the facets in $\partial F_J$, so we get a
list of the 1- and 0-faces along the way. The 3-faces are easily
deduced from the 2-faces.

The basis for our analysis is the following, which follows from the
theory of Gr\"obner bases (see~\cite{cohen} for instance, and also
section~\ref{sec:rur} of the present paper). Let $\ell$ be a number
field.
\begin{itemize}
  \item There is an algorithm to determine whether a system of $n$
    polynomial equations defined over $\ell$ in $n$ unknowns is
    0-dimensional (i.e. whether there are only finitely many solutions
    in $\C^n$);
  \item If the system is indeed 0-dimensional, there is an algorithm
    to determine the list of solutions; their entries lie in a finite
    extension $\k\supset \ell$. One can also determine the list of
    rational/real solutions.
  \item Polynomials with coefficients in $\ell$ can be evaluated at the
    solutions of a point with coordinates in $\k$, and one can
    determine whether the value is positive (resp. negative or zero).
\end{itemize}
When such systems have solution sets with unexpectedly high dimension,
there is usually a geometric explanation (typically some of the
intersecting bisectors share a slice, see~\cite{dpp2} for instance).
We will not address this issue, since it never occurs in the situation
of the present paper.

In all situations we will consider here, the extension $\ell$ will a
quadratic number field, and $\k$ will have degree at most four over
$\ell$. This makes all computations very quick (using capabilities of
recent computers, and standard implementations of Gr\"obner bases).

For the rest of the discussion, we make the following assumptions.
\begin{assume}\label{assumption1-2}
\begin{enumerate}
  \item For every $L\subset I$ with $|L|=4$, $F_L$ has dimension zero.
  \item For every $J\subset I$ with $|J|=2$, and every $x\in I$,
    $x\notin J$, the restriction $g_x$ of $f_x$ to $F_{J,\emptyset}$ has
    non-degenerate critical points.
\end{enumerate}
\end{assume}
These assumptions are by no means necessary in order to determine the
combinatorial structure of $F_{J,K}$, but they will simplify the
discussion in several places. Note also that they can be checked
efficiently using a computer, in particular we state
\begin{prop}
  Let $M$ be the figure eight knot complement. Then the Ford domain of
  the irreducible boundary unipotent representation
  $\rho:\pi_1(M)\rightarrow PU(2,1)$, centered at the fixed point of
  the holonomy of any peripheral subgroup satisfies
  assumption~\ref{assumption1-2}.
\end{prop}
In contrast, the domains that appear in~\cite{dpp2} do not satisfy
these hypotheses.

The combinatorial description of $F_J$ (i.e. its connected components,
and the list of facets adjacent to it) can be obtained by starting
from a description of $F_{J,\emptyset}$, and repeatedly studying
$F_{J,K\cup\{x\}}$ from $F_{J,K}$, where $x\in I$ is not in $J\cup K$.
The latter inductive step is done as follows.

The boundary $\partial F_{J,K}$ can be described as a union of arcs
contained in $F_{J\cup \{k\}, K\setminus \{k\}}$ for some $k\in
K$. For computational purposes, we will always assume that an arc is
homeomorphic to a closed interval, that its endpoints are vertices,
but none of its interior points are vertices. 

Note also that the arcs may not be equal to $F_{J\cup \{k\},
  K\setminus \{k\}}$, since $F_{J\cup\{x\},\emptyset}$ may have a
double point.

For each arc $a$ in $\partial F_{J,K}$ as above, we study the set
$$
F_{J\cup\{k,x\}},
$$ 
which by assumption~\ref{assumption1-2}(1) is obtained by solving a
0-dimensional system. Keeping only solutions that lie in $a$, we get a
subdivision of $a$ into connected components of $a\setminus F_{J\cup
  \{k,x\}}$, and for each such component we check whether or not it is
in $F_{J\cup\{k\},\{x\}}$. If so, it is a component of the boundary of
$F_{J,K\cup\{x\}}$.

We then compute the critical points of the restriction to $F_J$ of the
equation $f_x$ (this can be done because of
assumption~\ref{assumption1-2}(2), and determine whether any such
critical point is inside $F_{J,K}$.

Suppose $c$ is in a component $C_{J,K}$ of $F_{J,K}$. 
\begin{itemize}
\item If $g_x(c)=0$ and $c$ is a saddle point for the restriction
  $g_x$ of $f_x$, then a neighborhood of $c$ in
  $\overline{F}_{J,K\cup\{x\}}$ is the union of two sectors meeting in
  their apex. In a neighborhood of $c$, $F_{J,K\cup\{x\}}$ will have
  four boundary arcs. Each such arc will either connect $c$ to another
  saddle point of $g_x$, or it will connect it to a vertex in the
  boundary of $C_{J,K}$. For each such arc, we take a sample point to
  check whether it is contained in $F_{J\cup\{x\},K}$.
\item If $g_x(c)\neq 0$, there could be an isolated component of
  $F_{J\cup\{x\},K}$ that winds around $c$. In order to determine
  whether this happens or not, we consider the slice $t_1=\alpha_1$,
  and intersect it with $g_x=0$.  Recall that this intersection
  contains either 0, 1, or 2 points (because it is obtained by solving
  an equation that has degree at most two, which is not identically
  zero because $g_x(c)\neq 0$).

  Then there is an isolated component if and only if the intersection
  consists of precisely two points, and the two intersection points
  lie in the same connected component of $F_{J,K}$.
\end{itemize}
Now collecting the boundary arcs with the inside arcs (joining two
points that are either saddle or boundary vertices $F_{J\cap
  \{k,x\}}$), we get a stratum decomposition for
$F_{J,K\cup\{x\}}$. 

Moreover, if we make the following assumption, then all components of
$F_{J,K\cup\{x\}}$ are topological disks, since their boundary
consists of a single component.\\
\begin{assume}\label{assumption3}
\begin{enumerate}
  \setcounter{enumi}{2}
\item The curves $F_{J\cup\{x\},K}$ do not have any isolated component
  in $F_{J,K}$.
\end{enumerate}
\end{assume}
Once again, in the special case of the Ford domain relevant to the
irreducible boundary unipotent rank one, it turns out this hypothesis
is satisfied.


\subsubsection{Rational Univariate Representation} \label{sec:rur}

We briefly recall what we need about rational univariate
representations; for details on this technique, see~\cite{rouillier}.
Recall that given a 0-dimensional polynomial system
\begin{equation}\label{eq:sys1}
\left\{\begin{array}{l}
f(t_1,t_2)=0\\
g(t_1,t_2)=0
\end{array}\right.
\end{equation}
with coefficients in the number field $\ell$, we can write it as a
polynomial system with rational coefficients by using a primitive
element for $\ell$; the corresponding system has one more variable
(which we denote by $s$), and one more equation (which is the minimal
polynomial of a primitive generator for $\ell$). We write it in the
form
\begin{equation}\label{eq:sys2}
\left\{\begin{array}{l}
\tilde{f}(t_1,t_2,s)=0\\
\tilde{g}(t_1,t_2,s)=0\\
m(s)=0
\end{array}\right.
\end{equation}
where $\tilde{f}$ is obtained from $f$ by expressing its coefficients
as polynomials in the primitive element for $\ell$. In the cases that
interest us, $\ell$ will be a totally real number field, which we
assume from now on.

In this discussion we consider systems of two equations in two
variables (so we get 3 equations in 3 variables, counting the
extra-variable corresponding to the primitive element of the number
field), but we could also allow system that have more equations than
the number of variables (the important point is that the ideal
generated by the equations should be 0-dimensional).

Now the key point is that there exists a 1-variable polynomial $r$
such that the solutions are parametrized as rational functions of the
roots of $r$. More specifically, there exist polynomials $r$, $p_0$,
$p_1$, $p_2$ and $q$ with integer coefficients such that the solutions
of the system can be written in the form
\begin{equation}\label{eq:ratparam}
s=p_0(u)/q(u), t_1=p_1(u)/q(u),\quad t_2=p_2(u)/q(u),
\end{equation}
and the latter formula gives a solution of~\eqref{eq:sys2} if and
only $u$ is a root of $r$. Of course, since in general the minimal
polynomial $m$ has several roots, this produces more solutions of
system~\eqref{eq:sys1} than we would like. The solutions
of~\eqref{eq:sys1} can easily be obtained by sifting the solutions
of~\eqref{eq:sys2} once we know isolating intervals for the roots of
$m$.

Note that, even though all the equations relevant to this paper have
coefficients in a fixed number field (namely $\ell=\Q(\sqrt{7})$), the
vertices usually have entries in a larger number field (namely the
field generated by a given root of the rational parametrizing
polynomial $r$).

Note also that the solutions lie in a subfield $L\subset \mathbb{C}$
if and only if the corresponding root $u$ of $r$ lies in $L$. In
particular, if we want to find \emph{real} solutions of the system, we
can restrict to studying \emph{real} roots of $r$, which can be
specified by isolating intervals.

Using a rational univariate representation for the vertices provides a
convenient set of methods that allow us to:
\begin{enumerate}[(i)]
  \item find the list of faces that contain a given vertex;
  \item for each bounding bisector not containing a vertex, check
    which side the vertex is in;
  \item check if two vertices are the same;
  \item check whether a given vertex is inside a given arc;
  \item if two vertices in $F_{J\cup \{x\},\emptyset}$ are given,
    check whether these two vertices are joined by an arc in $F_{J\cup
      \{x\},\emptyset}$.
\end{enumerate}
Items~(i) and~(ii) are very simple because all our equations are
defined over a given $\ell$. 
Given a polynomial $h(t_1,t_2)=\tilde{h}(t_1,t_2,s)$, we start by
substituting the parametrization~\eqref{eq:ratparam} in $\tilde{h}$,
replacing $u$ by the appropriate interval of values of the rational
parameter. If the corresponding interval does not contain 0, we know
the sign of $h$ at that vertex. 

Otherwise, we keep the exact parametrization~\eqref{eq:ratparam} and get
a rational function in $u$ that represents $h$ at the solutions
of~\eqref{eq:sys2}, and we check whether it vanishes at the appropriate
root of $r$. This corresponds to checking whether our favorite root of
the rational parametrizing polynomial $r$ is also a root of another
given polynomial with integer coefficients (namely the numerator of
the above rational function); this can be done by computing their
greatest common divisor, and isolating its real roots. 

If the rational function does not vanish, we compute a more precise
interval for the value of $\tilde{h}$, and refine precision untill the
interval does not contain 0. Of course, in all generality, this may
require such high precision that it would exhaust the system memory,
but this does not seem to happen for the verifications that appear in
this paper, at least for our implementation on standard modern
computers.

We now sketch how to implement item~(iii). Suppose we are given two
rational parametrizations
$$
\begin{array}{l}
s=p_0(u)/q(u), t_1=p_1(u)/q(u),\quad t_2=p_2(u)/q(u)\\
s=a_0(v)/b(v), t_1=a_1(v)/b(v),\quad t_2=a_2(v)/b(v),
\end{array}
$$ 
where $u$ (resp. $v$) is to be taken to be a specific root of
$r(u)$ (resp. $c(v)$). Equality corresponds to verifying whether
$p_1(u)b(v)-q(u)a_1(v)$ (resp. $p_2(u)b(v)-q(u)a_2(v)$) vanishes at
the corresponding roots. If the rational parameters were the same,
this would simply amount to computing a greatest common divisor, but in
general the parameters from both rational representaions are different.

One way to handle this is to solve the system
$$
\left\{
\begin{array}{l}
p_1(u)b(v)-q(u)a_1(v)=0\\
r(u)=0\\
c(v)=0
\end{array}
\right.,
$$ 
which can be done using a rational univariate representation once
again. The result then follows from sifting solutions and keeping only
those that give the right root for $u$ and $v$, and checking whether
the sift gives a solution of not.

In order to explain how to check~(iv), we need to describe in more
detail how we encode arcs. We will assume
\begin{itemize}
\item that every arc is parametrized by one of the spinal coordinates
  (this can always be achieved, perhaps after subdividing certain arcs
  if necessary),
\item that the endpoints of every arc are vertices (parametrized by a
  rational univariate representation, as discussed above), and
\item that there are no vertices stricly inside any arc.
\end{itemize}

Then, in order to check whether a given vertex is inside an arc
parametrized by $t_1$, we need to compare its $t_1$ value with the
$t_1$ values of the endpoints of the arc. This amounts to checking the
sign of an expression of the form
$$
p_1(u)/q(u)-a_1(v)/b(v),
$$ 
where $u$ (resp. $v$) is a specific root of $r$ (resp. $c$). This
is the same as the test that occurs in item~(iii).

If the vertex $t_1$ value is between the $t_1$-values of the endpoints
of the arc, we still need to check whether it is in the correct
arc.

\subsubsection{Sample computations} \label{sec:sample}

We determine some sets $F_J$, $|J|=2$ explicitly, in order to
illustrate the phenomena that can occur when applying the algorithm
from the previous section. The general scheme to parametrize
$F_{J,\emptyset}$ is explained in~\cite{derauxfalbel}, for instance.

When $0\notin J=\{j,k\}$, we distinguish two basic cases, depending on
whether $p_0$, $p_j$ and $p_k$ are in a common complex line. This
happens if and only if some/any lifts $\tilde{p}_j\in\C^3$ are
linearly dependent. In that case, the bisectors $F_{\{j\}}$ and
$F_{\{k\}}$ have the same complex spine, and their intersection is
either empty or a complex line (this never happens in the Ford domains
studied in this paper).

Otherwise, $F_{J,\emptyset}$ can be parametrized by vectors of the form
$$
(\bar z_1 p_0 - p_j)\boxtimes (\bar z_2 p_0 - p_k) = z_1 p_{k0} + z_2 p_{0j} + p_{jk},
$$ 
with $|z_1|=|z_2|=1$, and where $p_{mn}$ denotes $p_m\boxtimes p_n$
(see section~\ref{sec:bisectors}).

Valid pairs $(z_1,z_2)$ in the Clifford torus $|z_1|=|z_2|=1$ are given by pairs where 
$$
\langle z_1 p_{k0} + z_2 p_{0j} + p_{jk}, z_1 p_{k0} + z_2 p_{0j} + p_{jk}\rangle <0,
$$
which can be rewritten as
$$
\Re(\mu_0(z_1)z_2)=\nu_0(z_1),
$$
for $\mu_0$ and $\nu_0$ affine in $z_1,\bar z_1$.

In terms of the notations of section~\ref{sec:algo}, the restriction
$g_0$ of $f_0$ to $F_{J,\emptyset}$ is given by
$$
g_0(z_1,z_2)=\Re(\mu_0(z_1)z_2)-\nu_0(z_1)).
$$

In order to draw pictures, we will sometimes use log-coordinates
$(t_1,t_2)$ for $F_{J,\emptyset}$, and write, for $j=1,2$,
$$
z_j=exp(2\pi i t_j).
$$

Given $l\notin J$, we already mentioned in section~\ref{sec:bisectors}
how to write the restriction $g_l$ of $f_l$ to to $F_J$. Note that $\langle
p_{k0},p_0\rangle=\langle p_{0j},p_0\rangle=0$, so the equation $f_x$
reads
$$
|\langle p_{jk}, p_0 \rangle| = |\langle z_1 p_{k0} + z_2 p_{0j} + p_{jk}, p_l \rangle|,
$$
which again can be written in the form
$$
\Re(\mu_l(z_1)z_2)=\nu_l(z_1).
$$

In order to compute the critical points of the restriction to
$|z_1|=|z_2|=1$ of a function $h(z_1,\bar z_1,z_2,\bar z_2)$, we
search for points where
$$
\left\{\begin{array}{l}
\frac{\partial h}{\partial z_1}z_1 - \frac{\partial h}{\partial \bar z_1}\bar z_1=0\\
\frac{\partial h}{\partial z_2}z_2 - \frac{\partial h}{\partial \bar z_2}\bar z_2=0
\end{array}\right.,
$$ 
Gr\"obner bases for the corresponding systems tell us whether these
critical points are non-degenerate (see assumption~\ref{assumption3}), and if so, we
can compute them fairly explicitly, i.e. describe their coordinates as
roots of explicit polynomials (in particular they can be computed to
arbitrary precision).

\begin{prop}\label{prop:12}
  Let $J=\{1,2\}$. Then $F_J$ is empty, and $\overline{F}_J$ is a
  singleton, given by $F_{\{1,2,3,5,10,11\}}$.
\end{prop}
The singleton in the Proposition is $\{p_2\}$, for $p_2$ as in
Lemma~\ref{lem:p2}. It follows from the Proposition that $p_2$ lies
precisely on six bounding bisectors (Lemma~\ref{lem:p2} only showed
that it was on at least six, listed in Tables~\ref{tab:vertices2}
and~\ref{tab:vertices3}).

\begin{pf}
For $J=\{1,2\}$, we get
$$
\mu_0(z_1)=-2-\bar z_1,\quad \nu_0(z_1)=-3+z_1+\bar z_1.
$$
The discriminant 
$$
|\mu|^2-\nu^2 = -6 + 16\Re z_1 - 2\Re z_1^2
$$
vanishes precisely for four complex values of $z_1$, which are the roots of 
\begin{equation}\label{eq:somepoly}
z_1^4 - 8 z_1^3 + 6 z_1^2 - 8 z_1 + 1. 
\end{equation}
Since we know $F_{J,\{0\}}$ is connected (see~\cite{goldman},
Theorem~9.2.6), we know that at most two of these roots lie on the
unit circle. In fact, $z_1=z_2=1$ gives a point in $F_{J,\{0\}}$,
so $F_{J,\{0\}}$ is non empty, hence there must be two (complex
conjugate) roots on the unit circle. Indeed, these roots have argument
$2\pi t$ with $t=\pm 0.20682703\dots$.

A more satisfactory way to check that the
polynomial~\eqref{eq:somepoly} has precisely two roots on the unit
cirle is to split $z_1=x_1+iy_1$ into its real and imaginary parts
(this gives a general method that does not rely on geometric
arguments).

Indeed $z_1$ is a root of~\eqref{eq:somepoly} if and only if
$(x_1,y_1)$ is a solution of the system $-6+16x_1-2x_1^2+2y_1^2=0$,
$x_1^2+y_1^2=1$. These equations imply that $x_1=2\pm\sqrt{3}$, and then
$$
y_1^2=2-4x_1,
$$ 
which is positive only for $x_1=2-\sqrt{3}$, and then we get
$y_1=\pm\sqrt{4\sqrt{3}-6}$.

In order to run the algorithm from the preceding section, we write the
restriction $g_3$ of $f_3$ to $F_{J,\emptyset}$, which is given by 
$$ 
-3 + 2\Re\{ \frac{1-i\sqrt{7}}{2} z_1 + \frac{5-i\sqrt{7}}{2} z_2 +
\frac{-3+i\sqrt{7}}{2} z_1\bar z_2 \}.
$$
Gr\"obner bases calculations show that the system
$g_0(z)=g_3(z)=|z_1|^2-1=|z_2|^2-1=0$ has precisely two solutions,
given in log-coordinates by
$$
(-0.20418699\dots,-0.03294828\dots),\quad (0.15576880\dots,-0.07655953\dots).
$$ 
Once again, the most convenient way to use Gr\"obner bases is to
work with four variable $x_1,y_1,x_2,y_2$ given by real and imaginary
parts of $z_1$ and $z_2$ (with extra equations $x_j^2+y_j^2=1$.

The combinatorics of $F_{J,K}$ for $K=\{0,3\}$ are illustrated in
Figure~\ref{fig:algo-12}(b). It is a disk with two boundary arcs, given by
$F_{\{1,2,0\},\{3\}}$ and  $F_{\{1,2,3\},\{0\}}$.

\begin{figure}
  \subfigure[$F_{\{1,2\},\{0\}}$]{\epsfig{figure=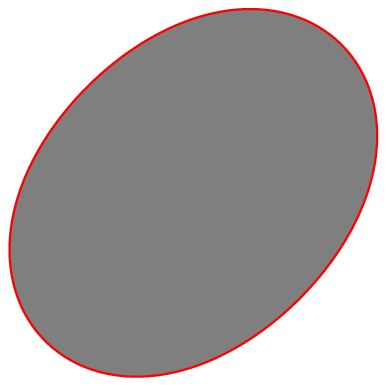,height=0.15\textwidth}}\hfill
  \subfigure[$F_{\{1,2\},\{0,3\}}$]{\epsfig{figure=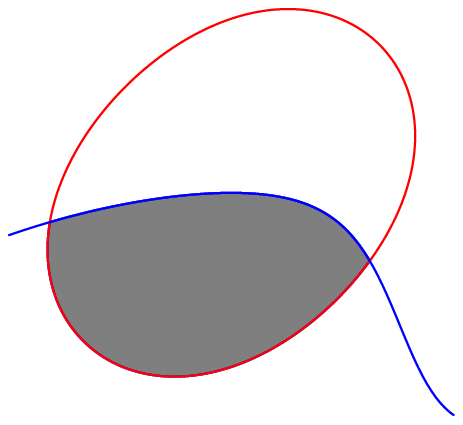,height=0.15\textwidth}}\hfill
  \subfigure[$F_{\{1,2\},\{0,3,5\}}$]{\epsfig{figure=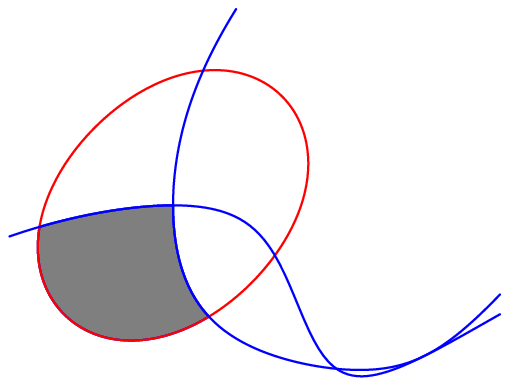,height=0.15\textwidth}}\hfill
  \subfigure[$F_{\{1,2\}}$]{\epsfig{figure=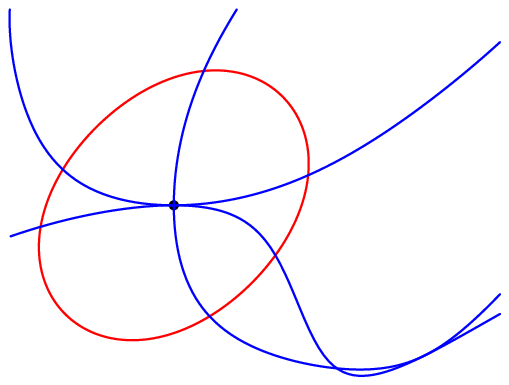,height=0.15\textwidth}}
  \caption{Steps of the algorithm to determine $F_{\{1,2\}}$.}\label{fig:algo-12}
\end{figure}

As the next element to include in $K$, we choose 5 rather than 4, in
order to shorten the discussion slightly. The curve $F_{\{1,2,5\},\emptyset}$
intersects $F_{\{1,2,0\},\emptyset}$ two points, given in log-coordinates by
$$
(0.04600543\dots, 0.20593006\dots),\quad (0.05483483\dots, -0.17019919\dots).
$$
Only the second one is inside the arc $F_{\{1,2,0\},\{3\}}$.

The curve $F_{\{1,2,5\},\emptyset}$ intersects
$F_{\{1,2,3\},\emptyset}$ in five points, given by $(z_1,z_2)=$
$$
(1,1),\ (i,-i),\ (-i,i), (\frac{9+5i\sqrt{7}}{16},\frac{-3+i\sqrt{7}}{4}),\ 
(\frac{-3+i\sqrt{7}}{4},\frac{1-3i\sqrt{7}}{8}).
$$
only one of which is in $F_{\{1,2,3\},\{0\}}$, namely $(1,1)$.

Now $F_{\{1,2\},\{0,3,5\}}$ has three boundary arcs, given by
  $F_{\{1,2,0\},\{3,5\}}$, $F_{\{1,2,3\},\{0,5\}}$ and
      $F_{\{1,2,5\},\{0,3\}}$ (see Figure~\ref{fig:algo-12}(c)).

Next, we include 10 in $K$. The curve $F_{\{1,2,10\},\emptyset}$ intersects
$F_{\{1,2,0\},\emptyset}$ in two points, none of which is in
$F_{\{1,2,0\},\{3,5\}}$. Hence the arc $F_{\{1,2,0\},\{3,5\}}$ is either
completely inside, or completely outside
$F_{\{1,2,0\},\{3,5,10\}}$. One easily checks that it is outside, by
taking a sample point.

The curve $F_{\{1,2,10\},\emptyset}$ intersects
$F_{\{1,2,3\},\emptyset}$ in five points, none of which is in
$F_{\{1,2,3\},\{0,5\}}$. The arc $F_{\{1,2,3\},\{0,5\}}$ is either
completely inside, or completely outside $F_{\{1,2,3\},\{0,5,10\}}$
and a sample point shows it is outside.

Similarly, the curve $F_{\{1,2,10\},\emptyset}$ intersects
$F_{\{1,2,5\},\emptyset}$ in six points, none of which is in
$F_{\{1,2,5\},\{0,3\}}$, and the arc $F_{\{1,2,5\},\{0,3\}}$ is
completely outside $F_{\{1,2,3\},\{0,5,10\}}$.

This implies that $F_{\{1,2\}}$ is empty (see Figure~\ref{fig:algo-12}(d)).

Finally we consider the intersection of $F_{\{1,2,10\},\emptyset}$
with the three vertices of $F_{\{1,2\},\{0,3,5\}}$. One easily checks
that the only intersection is the point with complex spinal
coordinates given by $(1,1)$, and this point indeed a vertex of $F$.
It is in homogeneous coordinates in $\C^3$ given by
$$
(\frac{3-i\sqrt{7}}{2},-2,-\frac{3-i\sqrt{7}}{2}),
$$ 
and that it is on precisely six bounding bisectors (by construction
it is on $\B_1$ and $\B_2$, and it is also in $\B_3$, $\B_5$,
$\B_{10}$ and $\B_{11}$).
In terms of the notation of section~\ref{sec:algo}, this point is
$$
F_{\{1,2,3,5,10,11\}}.
$$ 
In fact one easily checks that this point is the fixed point of $G_2$
(which by definition of the bounding bisectors is obviously in
$\B_1\cap\B_2$).
\end{pf}

\begin{rmk}
  \begin{enumerate}
    \item Throughout the proof of Proposition~\ref{prop:12}, we have
      ignored the issue of critical points. In principle, at each
      stage, we may have missed some isolated components of the curves
      $F_{\{1,2,k\},\emptyset}$; if this were the case, the set
      $F_{\{1,2\}}$ would still be contained in the set which we just
      described, hence it must be empty anyway.
    \item The curves $F_{\{1,2,10\},\emptyset}$ and
      $F_{\{1,2,3\},\emptyset}$ are in fact tangent at $(1,1)$, which
      is a vertex of $F$. We shall come back to this point later, when
      discussing stability of the combinatorics of $F$ under
      deformations.
  \end{enumerate}
\end{rmk}

\begin{prop} \label{prop:13}
  $F_{\{1,3\}}$ is combinatorially a triangle, with three boundary
  arcs given by $F_{\{1,3,0\}}$, $F_{\{1,3,5\}}$, $F_{\{1,3,11\}}$,
  and three vertices given by $F_{\{0,1,3,5\}}$, $F_{\{0,1,3,11\}}$,
  and $F_{\{1,2,3,5,10,11\}}$.
\end{prop}

Note that this triangle appears in Figure~\ref{fig:face2}(a)
and~\ref{fig:face3}(a), it is the intersection of the bounding
bisectors $\B_1$ (resp. $\B_3$) corresponding to $G_2$
(resp. $G_3$). The edges in $\ch 2$ are on $\B_5$, which corresponds
to $G_1G_2G_1^{-1}$ and $\B_{11}$, which correspond to
$G_1^{-1}G_3G_1$.

\begin{pf}
  As in the argument for $F_{\{1,2\}}$, we study $F_{J,K}$ for
  increasing sets $K$, freely choosing the order we use to increase
  $K$. We describe an efficient way to get down to
  $F_{\{1,3\}}$ in the form of a picture, see Figure~\ref{fig:algo-13}.
  \begin{figure}[htbp]
    \subfigure[$F_{\{1,3\},\{0,5\}}$]{\epsfig{figure=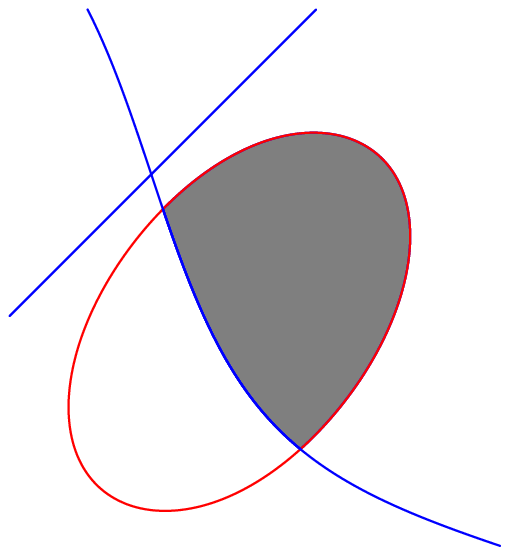,height=0.3\textwidth}}\hfill
    \subfigure[$F_{\{1,3\},\{0,5,11\}}$]{\epsfig{figure=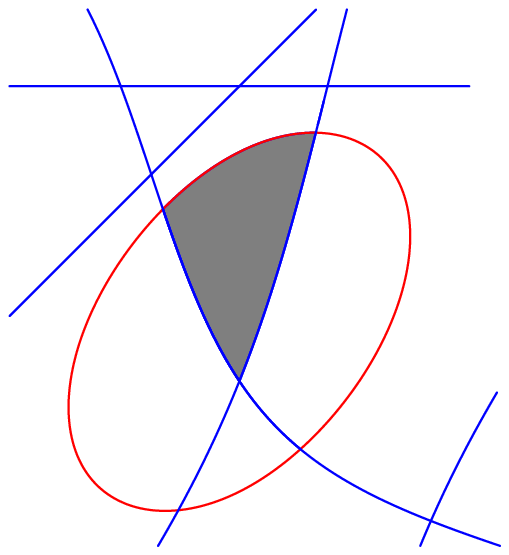,height=0.3\textwidth}}\hfill
    \subfigure[$F_{\{1,3\},\{0,5,11,2,10\}}$]{\epsfig{figure=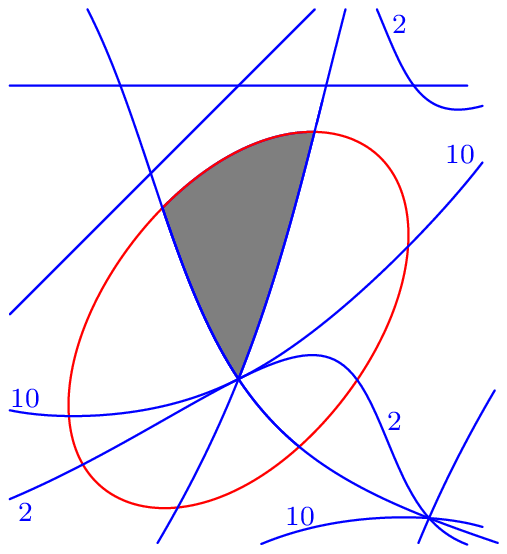,height=0.3\textwidth}}
    \caption{Steps of the algorithm to determine $F_{\{1,3\}}$.}\label{fig:algo-13}
  \end{figure}

  We start by studying $F_{\{1,3\},\{5\}}$. Note that the curve
  $F_{\{1,3,5\},\emptyset}$ has two double points. These points can be
  obtained by writing the equation $g_5$ as
  $$
  \Re(\mu(z_1)z_2) = \nu(z_1),
  $$
  where 
  $$
  \mu(z_1)=\frac{3+i\sqrt{7}}{2}-\bar z_1,\quad \nu(z_1)=1-\Re(\frac{3+i\sqrt{7}}{2}z_1).
  $$
  The discriminant $|\mu(z_1)|^2-\nu(z_1)^2$ is given by
  $$
  2 + \Re\left( \frac{-1-3i\sqrt{7}}{4}z_1^2  \right),
  $$
  which vanishes for 
  $z_1=\pm\frac{3-i\sqrt{7}}{4}.$
  Plugging this back into the equation $g_5$ gives 
  $z_2=\mp\frac{3-i\sqrt{7}}{4}.$
  One easily checks that $g_0(z_1,z_2)>0$ for these two double
  points, i.e. they lie outside complex hyperbolic space.
  
  One checks that $F_{\{1,3,5\},\emptyset}$ intersects
  $F_{\{1,3,0\},\emptyset}$ in precisely two points (and these
  intersections are transverse), so we get two arcs in the boundary of
  $F_{\{1,3\},\{0,5\}}$, namely $F_{\{1,3,5\},\{0\}}$ and
  $F_{\{1,3,0\},\{5\}}$ (see Figure~\ref{fig:algo-13}(a)). 

  In principle, there could be an extra arc in $F_{\{1,3,5\},\{0\}}$,
  not intersecting $F_{\{1,3,0\},\emptyset}$, so we compute critical
  points of $g_5$. Their are given by solutions of the system
  $$
  \left\{\begin{array}{l}
  \Im\{ (\bar z_2 +\frac{3+i\sqrt{7}}{2})z_1 \} = 0\\
  \Im\{ (\bar z_1 +\frac{3+i\sqrt{7}}{2})z_2 \} = 0
  \end{array}\right.,
  $$
  that satisfy $|z_1|=|z_2|=1$.
  
  There are four such critical points, of the form $(\pm
  \alpha,\pm\alpha)$ where $\alpha=\frac{3-i\sqrt{7}}{4}$ (of course
  this list includes the double points computed before). The
  corresponding points are outside $F$, in fact
  $g_0(\pm\alpha,\pm\alpha)>0$.

  A similar analysis justifies part~(b) of Figure~\ref{fig:algo-13},
  i.e. that $\overline{F}_{\{1,3\},\{0,5,11\}}$ is combinatorially a
  triangle (with one side on $\partial_\infty \ch 2$).

  We sketch how to justify that
  $F_{\{1,3\}}=F_{\{1,3\},\{0,5,11\}}$. For $k=2$ and $k=10$, the
  curve $F_{\{1,3,k\},\emptyset}$ actually goes through a vertex of
  $F_{\{1,3\}}=F_{\{1,3\},\{0,5,11\}}$; for $k\neq 0,2,5,10,11$,
  $F_{\{1,3,k\},\emptyset}$ does not intersect even
  $\overline{F}_{\{1,3\},\{0,5,11\}}$. 

  We start by studying $F_{\{1,3,0\},\emptyset}\cap
  F_{\{1,3,2\},\emptyset}$. In order to use standard root isolation
  methods, we use real equations, in $x_1,y_1,x_2,y_2$. Computing a
  Gr\"obner basis for the ideal generated by the equations $g_0$,
  $g_3$, $x_1^2+y_1^2-1$ and $x_2^2+y_2^2-1$, we see that it contains
  $$
  39-840\sqrt{7}y_2+4088y_2^2+608y_2^3\sqrt{7}-9152y_2^4+1024y_2^5\sqrt{7}+7168y_2^6,
  $$
  which has precisely two real roots, given approximately by
  $y_2^{(1)}=0.01815877\dots$ and $y_2^{(2)}=0.65602473\dots$. 

  The Gr\"obner basis also gives an expression for $x_1,y_1,x_2$ in terms of
  $y_2$, namely 
  $$
  x_1=\{ -4943+16836\sqrt{7}y_2-142640y_2^2+53184y_2^3\sqrt{7}+72128y_2^4-75264y_2^5\sqrt{7} \}/14725,
  $$
  $$
  y_1=\{ 5058\sqrt{7}+45888y_2-112560y_2^2\sqrt{7}+309472y_2^3+74432y_2^4\sqrt{7}-422912y_2^5 \}/14725,
  $$
  $$
  x_2=\{ 20-21\sqrt{7}y_2+16y_2^2+32y_2^3\sqrt{7} \}/19.
  $$

  Substituting either value $y_2^{(j)}$ gives two points
  $a^{(j)}=(x_1^{(j)},y_1^{(j)},x_2^{(j)},y_2^{(j)})$, $j=1,2$ and we
  claim that $g_5(a^{(1)})>0$ and $g_{11}(a^{(2)})>0$. Clearly this
  can be checked by simple interval arithmetic, in fact
  $$
  g_5(a^{(1)})=3.80716606\dots, g_{11}(a^{(2)})=3.94518313\dots.
  $$

  The analysis of $F_{\{1,3,5\},\emptyset}\cap
  F_{\{1,3,2\},\emptyset}$ is in a sense simpler, because all the
  solutions to the corresponding system are defined over
  $\mathbb{Q}(i,\sqrt{7})$. The system has precisely five solutions, given by
  \begin{equation*}
    \begin{array}{c}
      (i,\frac{1+\sqrt{7}}{4}+i\frac{1-\sqrt{7}}{4}), 
      (-i,\frac{1-\sqrt{7}}{4}-i\frac{1+\sqrt{7}}{4}),\\ 
      (\frac{-3+i\sqrt{7}}{4},\frac{3-i\sqrt{7}}{4}), 
      (\frac{9+5i\sqrt{7}}{16},-\frac{9+5i\sqrt{7}}{16},
      (1,\frac{3+i\sqrt{7}}{4}).
    \end{array}
  \end{equation*}
  Only one of these solutions satisfies $g_0\leq 0$, namely the
  last one (in other words, only one intersection point lies $\chb
  2$).
  
  Note that we already found one point in $F_{\{1,3,2\},\emptyset}\cap
  F_{\{1,3,5\},\emptyset}$, namely the fixed point of $G_2$ (see the
  proof of Proposition~\ref{prop:12}). 

  Similarly, one verifies that $F_{\{1,3,2\},\emptyset}\cap
  F_{\{1,3,11\},\emptyset}$ contains precisely six points, only one of
  which gives a point in (the closure of) complex hyperbolic
  space.

  Once again, since we already know one point in this intersection
  (namely the fixed point of $G_2$), we get that
  $F_{\{1,3,2\},\emptyset}$ with $\partial
  F_{\{0,1,3,5,11\},\emptyset}$ consists of precisely one point. This
  implies that $\partial F_{\{0,1,3,5,11\},\emptyset}$ is either
  completely inside or completely outside $\partial
  F_{\{0,1,3,5,11,2\},\emptyset}$. It is easy to check that it is
  inside, by testing a sample point (for instance one of the other
  vertices of the triangle $\partial F_{\{0,1,3,5,11\},\emptyset}$).

  We now show that $F_{\{1,3,2\},\emptyset}$ does not intersect
  $F_{\{0,1,3,5,11\},\emptyset}$, by computing the critical points of
  $g_2$. There are six critical points, given by
  \begin{equation*}
    \begin{array}{c}
      (-1,-\frac{1+3i\sqrt{7}}{8}), 
      (\frac{3-i\sqrt{7}}{4},\frac{3-i\sqrt{7}}{4}),
      (\pm \frac{1+i\sqrt{7}}{\sqrt{8}},  \pm \frac{1-i\sqrt{7}}{\sqrt{8}}),
    \end{array}
  \end{equation*}  
  and one easily checks that none of them is inside
  $F_{\{0,1,3,5,11\},\emptyset}$. In particular, we get that the
  minimum value of $g_3$ on $\overline{F}_{\{0,1,3,5,11\},\emptyset}$
  is $0$, and it is realized precisely at one vertex (namely the fixed
  point of $G_2$).

  In other words, we get
  $F_{\{0,1,3,5,11\},\emptyset}=F_{\{0,1,2,3,5,11\},\emptyset}$,
  i.e. including the inequality $g_2<0$ at this stage has no effect.
  An entirely similar computation shows that
  $F_{\{0,1,2,3,5,11\},\emptyset}=F_{\{0,1,2,3,5,10,11\},\emptyset}$.

  For all $k\neq 0,1,2,3,5,10,11$, $F_{\{0,1,3,k\},\emptyset}$ does
  not intersect even the closure
  $\overline{F}_{\{0,1,2,3,5,11\},\emptyset}$, and one can use
  arguments as above using interval arithmetic.
\end{pf}

Similar arguments allow us to handle the detailed study of all the
polygons that appear on Figure~\ref{fig:face2} and~\ref{fig:face3}.
\begin{prop}\label{prop:14}
  $F_{\{1,4\},\emptyset}$ is a Giraud disk, which is entirely
  contained in the exterior of $\B_5$. In particular, $F_{\{1,4\}}$ is
  empty.
\end{prop}

\begin{pf}
  We will prove that $F_{\{5\},\emptyset}$ does not intersect the
  Giraud torus $\widehat{F}_{\{1,4\},\emptyset}$. In order to see this, we use
  complex spinal coordinates, and write $g_5(z_1,z_2)$ for the
  restriction of $f_5$ to the Clifford torus $|z_1|=|z_2|=1$.

  One computes explicitly that 
  $$
  g_5(z_1,z_2)=4+2\Re\{\frac{1+i\sqrt{7}}{2}z_1\bar z_2\}.
  $$
  This is clearly always positive when $|z_1|=|z_2|=1$.

  In other words, the Giraud torus $\widehat{F}_{\{1,4\},\emptyset}$
  is entirely outside $F$.
\end{pf}

\begin{prop}\label{prop:16}
  $F_{\{1,6\},\emptyset}$ is empty. The Giraud torus
  $\widehat{F}_{\{1,6\},\emptyset}$ is completely outside complex
  hyperbolic space, in other words the bisectors $\B_1$ and $\B_6$ are
  disjoint.
\end{prop}

\begin{pf}
  We write the equation of $F_{\{0,1,6\},\emptyset}$ in spinal
  coordinates for the Giraud torus $F_{\{1,6\},\emptyset}$, which reads
  $$
  g_0(z_1,z_2)=18-2\Re\{ 4(z_1+z_2) + z_1\bar z_2\}.
  $$ 
  Clearly this is non-negative when $|z_1|=|z_2|=1$, and in that case
  it is zero if and only if $z_1=z_2=1$.
  
  In other words, $\widehat{\B}_1$ and $\widehat{\B}_2$ intersect in
  a point in $\chb 2$.
    Note that this point is not in the closure of $F$, in fact it is
  strictly outside the half spaces bounded by $\B_2$, $\B_3$, $\B_5$,
  $\B_7$ and $\B_{11}$.
\end{pf}

\begin{prop}\label{prop:38}
  $F_{\{3,8\}}$ is empty. The Giraud torus $F_{\{3,8\},\emptyset}$
  contains a disk in $\ch 2$, but $\overline{F}_{\{3,8\},\{2,6\}}$ is
  empty.
\end{prop}

\begin{pf}
  The proof is actually very similar to that of
  Proposition~\ref{prop:13}, but since the corresponding set is empty,
  we go through some of the details.

  The curve $F_{\{3,8,2\},\emptyset}$ intersects
  $F_{\{3,8,0\},\emptyset}$ in precisely two points, and cuts out a
  disk in the Giraud disk $F_{\{3,8\},\emptyset}$, so that
  $F_{\{3,8\},\{0,2\}}$ is a disk with only two boundary arcs.

  One then easily verifies that $F_{\{3,8,6\},\emptyset}$ does not
  intersect $\overline{F}_{\{3,8\},\{0,2\}}$, so
  $F_{\{3,8\},\{0,2,6\}}$ is either equal to $F_{\{3,8\},\{0,2\}}$ or
  is empty (one needs to check critical points in order to verify
  this). 

  By taking a sample point $z$, and checking $f_6(z)>0$, one gets that
  $F_{\{3,8\},\{0,2,6\}}$ is empty.
\end{pf}

The study of $\B_1\cap\B_k$ for various values of $k$ is similar to
one of the previous few propositions, we list the relevant arguments
in Table~\ref{tab:face2}. When the proof is similar to
Proposition~\ref{prop:14}, the indices $l$ listed in brackets indicate
that $\B_1\cap\B_k$ is entirely outside the half space bounded by
$\B_l$.

The corresponding list of arguments used to study of $\B_3\cap\B_k$
for various values of $k$ in Table~\ref{tab:face3}.

Note that the arguments for $\B_2$ (resp. $\B_4$) are of course almost
the same as those for $\B_1$ (resp. $\B_3$), since the corresponding
faces are actually paired by $G_2$ (resp. $G_3$).

\begin{table}[htbp]
  \begin{tabular}{|c|c|}
    \hline 
    Prop~\ref{prop:finite} & 8, 14-16, 21-25, 29-33, 35\\ 
    Prop~\ref{prop:12} & 2, 12, 19, 26\\ 
    Prop~\ref{prop:13} & 3, 5, 9, 10, 11, 18, 20, 28\\ 
    Prop~\ref{prop:14} & 4[5,10], 7[3], 13[2,5,10], 17[9], 27[9,18], 36[17,28,34]\\ 
    Prop~\ref{prop:16} & 6, 34\\\hline
  \end{tabular}
  \caption{We list the indices where the arguments of each proposition
    apply to study $\B_1\cap\B_k$.}\label{tab:face2}
\end{table}

\begin{table}[htbp]
  \begin{tabular}{|c|c|}
    \hline
    Prop~\ref{prop:finite} & 16, 17, 22-36\\
    Prop~\ref{prop:12}     & 10, 13\\
    Prop~\ref{prop:13}     & 1, 2, 5, 6, 7, 11\\
    Prop~\ref{prop:14}     & 9[11], 14[7], 15[7], 18[1,10], 19[11], 20[7], 21[6,13]\\
    Prop~\ref{prop:16}     & 4, 12\\    
    Prop~\ref{prop:38}     & 8\\\hline
  \end{tabular}
  \caption{We list the indices where the arguments of each proposition
    apply to study $\B_3\cap\B_k$.}\label{tab:face3}
\end{table}

\subsubsection{Genericity}\label{sec:generic}

In order study deformations $\rho_t$ of the boundary unipotent
representation $\rho_0:\pi_1(M)\rightarrow PU(2,1)$, we will need more
information that just the combinatorics.

We will determine the non-transverse bisector intersections, and prove
that they remain non-transverse in the family of Ford domains for
groups in the 1-parameter family where the unipotent generator becomes
twist parabolic.

The basic fact is the following, which follows from the restrictive
character of the bounding bisector, namely they are all covertical
(because they define faces of a Ford domain).
\begin{prop}
  Let $J=\{j,k\}$ with $j\neq k$. Then the intersection
  $F_{\{j\},\emptyset}\cap F_{\{k\},\emptyset}=F_{J,\emptyset}$ is
  transverse at every point of $F_{J,\emptyset}$.
\end{prop}

The analogous statement is not true when $|J|\geq 3$, since
$F_{J,\emptyset}$ can have singular points (see
Figure~\ref{fig:algo-12} for instance). This will not be bothersome in
the context of our polyhedron $F$, because of the following.
\begin{prop}
  Suppose $|J|=3$ and $F_{J}$ is non-empty. Then the corresponding
  intersection of three bisectors (or two bisectors and
  $\partial_\infty \ch 2$) is transverse at every point of $F_J$.
\end{prop}
\begin{pf}
  This follows from the fact that double points of $F_{J,\emptyset}$
  occur only away from the face $F_J$. Indeed, one can easily locate
  these double points by the techniques explained in
  section~\ref{sec:sample}, and check that they are outside $F$ by
  using interval arithmetic.
\end{pf}

The situation near vertices is slightly more subtle, mainly because
our group contains some torsion elements, hence one expects the
intersections to be non-generic near the fixed points of those torsion
elements.

We will check possible tangencies between 1-faces intersecting at each
vertex. More generally, for each $j\neq k$, we will study tangencies
between all the curves of the form $F_{\{j,k,l\},\emptyset}$ for $\neq
j,k$ that occur at a vertex of $F$.

\begin{prop}
  Let $p$ be an ideal vertex of $F$, i.e. a vertex in $\partial_\infty
  \ch 2$. Then there are precisely three bounding bisectors $\B_i$,
  $\B_j$ and $\B_k$ meeting at $p$ (where $i,j,k>0$). The intersection
  of the four hypersurfaces in $\C^2$ given by the three extors
  $\widehat{\B_i}$, $\widehat{\B_j}$, $\widehat{\B_k}$, and
  $\partial_\infty \ch 2$ is transverse; in particular, none of the
  four incident 1-faces are tangent at $p$.
\end{prop}

Note that the ideal 1-faces are drawn in red on
Figures~\ref{fig:face2} and~\ref{fig:face3}, so the vertices on the
red curves are the ideal ones. The indices $(i,j,k)$ that appear in
the Proposition, i.e the bounding bisectors that contain a given ideal
vertex can be read off Figure~\ref{fig:hexagons}. For examples,
$(1,3,5)$, $(1,3,11)$, $(1,9,11)$,\dots are triples of indices that
correspond to ideal vertices.

\begin{pf}
  We treat the example of $F_{\{0,1,3,5\}}$, the other ones being
  entirely similar. The parametrization of the Giraud disk
  $F_{\{1,3\},\{0\}}$ was already explained in
  section~\ref{sec:sample}.

  The relevant vertex satisfies
  \begin{equation}\label{eq:idealvertex}
  x_1=0.80979557\dots, y_1 = -0.58671213\dots, x_2 = -0.53336432\dots, y_2 = 0.84588562\dots
  \end{equation}

  We write the equations of the bisectors in affine coordinates for
  complex hyperbolic space corresponding to the spinal coordinates,
  i.e. such that $(z_1,z_2)$ corresponds to
  $$
  p_{13} + z_1 p_{30} + z_2 p_{01},
  $$
  where $p_{jk}$ stands as before as the box product $p_j\boxtimes p_k$.

  In these coordinates, $\B_1$ is given by $|z_1|=x_1^2+y_1^2=1$, and $\B_3$ is
  given by $|z_2|=x_2^2+y_2^2=1$, and of course other bisectors have more
  complicated equations.

  The equation of the boundary of the ball is
  $$
  2 - \sqrt{7}y_2 - 4x_1 - x_2  - y_2\sqrt{7}x_1 + x_2\sqrt{7}y_1 + 2x_1^2 + 2y_1^2 + x_2^2 + y_2^2 -x_2x_1 -y_2y_1, 
  $$
  and the equation for $\B_5$ is given by
  $$
  3(x_1+x_2) - \sqrt{7}(y_1+y_2) - 2 x_2 x_1 - 2 y_2y_1 - x_1^2 - y_1^2  -x_2^2 -y_2^2.
  $$
  One then computes the gradient of each of these four equations, and
  checks that they are linearly independent at the point from
  equation~\eqref{eq:idealvertex} (this is readily done using interval
  arithmetic).
\end{pf}

\begin{prop}\label{prop:transversevertices}
  There are precisely six bounding bisectors containing $p_2$, indexed
  by 1,2,3,5,10,11. The pairwise and 3-fold intersections of these six
  bisectors are all transverse, but some 4-fold are not, namely
  \{1,2,3,10\}, \{1,2,5,11\}, \{3,5,10,11\}.
\end{prop}
  The precise list of bisectors that contain this vertex were already
  justified in section~\ref{sec:fixedpoints}, see Lemma~\ref{lem:p2}
  and Proposition~\ref{prop:12}. The point of
  Proposition~\ref{prop:transversevertices} is to give precise
  information about transversality. Recall from
  section~\ref{sec:fixedpoints} that $p_2$ is by definition the
  isolated fixed point of $G_2$, and the bisectors $\B_1$, $\B_2$,
  $\B_3$, $\B_5$, $\B_{10}$ and $\B_{11}$ are the bounding bisectors
  corresponding to the group elements $G_2$, $G_2^{-1}$, $G_3$,
  $G_1G_2$, $G_1^{-1}G_2^{-1}$, $G_1^{-1}G_3$, respectively (see
  section~\ref{sec:statecombo}).

\begin{pf}
  We work in spinal coordinates for $\B_1\cap \B_3$, and as in the
  preceding proof, we use $z_j=x_j+iy_j$, $j=1,2$ as global
  coordinates on $\ch 2$. The point $p_2$ is given by $z_1=1$,
  $z_2=\frac{3+i\sqrt{7}}{4}$.

  The equations of the six bisectors are as follows:
  $$
  \begin{array}{|r|l|}
    \hline
    1 & 4-4(x_1^2+y_1^2)\\
    2 & 2 + x_1 + 2x_2 + (y_1 - 2y_2)\sqrt{7} + (x_1y_2 - x_2y_1)\sqrt{7} + 3(x_1x_2 + y_1y_2)  - (x_1^2 + y_1^2) - 4(x_2^2+y_2^2)\\
    3 & 4-4(x_2^2+y_2^2)\\
    5 & 3(x_1 + x_2) - \sqrt{7}(y_1+y_2) - 2(x_2x_1 +y_2y_1) - (x_1^2 + y_1^2) - (x_2^2+y_2^2)\\
    10 & 2 - 4(x_1-x_2) + 4(x_2x_1 + y_2y_1) - 2(x_1^2 +y_1^2) - 2(x_2^2 +y_2^2) \\
    11 & 3 - 2x_2 + 3x_1 + \sqrt{7}y_1 + 3(x_1x_2 + y_1y_2)  + (x_2y_1 - y_2x_1)\sqrt{7} - 4(x_1^2+y_1^2) - (x_2^2 + y_2^2)\\
    \hline
  \end{array}
  $$
  One computes the gradients at the point $x_1=1$, $y_1=0$, $x_2=3/4$, $y_2=\sqrt{7}/4$, which are given by
  $$
  \begin{array}{c}
    v_1 = (-8, 0, 0, 0)\\
    v_2 = (3, \sqrt{7}, -1, -3\sqrt{7})\\
    v_3 = (0, 0, -6, -2\sqrt{7})\\
    v_5 = (-1/2, -3\sqrt{7}/2, -1/2, -3\sqrt{7}/2)\\
    v_{10} = (-5, \sqrt{7}, 5, -\sqrt{7})\\
    v_{11} = (-9/2, 5\sqrt{7}/2, -1/2, -3\sqrt{7}/2)
  \end{array},
  $$
  and the claim of the proposition follows from explicit rank computations.

  The tangent vectors to the intersection are given by
  $$
  \begin{array}{c}
    u_1 = (0,8/3,-\sqrt{7}/3,1)\\
    u_2 = (0,0,-3\sqrt{7},1)\\
    u_3 = (-2\sqrt{7}/3,-2/3,-\sqrt{7}/3,1)\\
  \end{array},
  $$ 
  and one easily checks that any curve tangent to these vectors must
  exit the polyhedron in a transverse fashion, more specifically, the
  exited bisectors are given in Table~\ref{tab:exits}.
  \begin{table}[htbp]
  $$
  \begin{array}{c|c|c|c}
    \textrm{vector} & \textrm{tangent to} & \textrm{exit in + direction} & \textrm{exit in - direction}\\\hline
    u_1 & 1,2,3,10      & 11        & 5\\
    u_2 & 1,2,5,11      & 3         & 10\\
    u_3 & 3,5,10,11     & 1         & 2
  \end{array}
  $$
  \caption{Each direction tangent vector $u_k$ to a non-tranverse
    quadruple intersection at $p_2$ exits the polyhedron; in the last
    two columns we list the two half spaces it exits (transversely) in
    the $\pm u_k$ direction.}\label{tab:exits}
  \end{table}
\end{pf}

\section{Side pairings}\label{sec:pairings}

\subsection{Faces paired by $G_2$}

We now justify the fact that $G_2^{-1}$ defines an isometry between
the faces for $G_2$ and $G_2^{-1}$. On the level of 2-faces, this
follows from the following.

\begin{prop}\label{prop:pairing2}
  The isometry $G_2^{-1}$ maps 
  \begin{enumerate}
    \item $G_3 p_0$ to $G_1^{-1}G_3 p_0$;
    \item $G_1^{-3}G_3^{-1} p_0$ to $G_1G_3^{-1} p_0$;
    \item $G_1^{-1}G_3 p_0$ to $G_1^{-1}G_2^{-1} p_0$;
    \item $G_1^{-1}G_2 p_0$ to $G_3^{-1}p_0$;
    \item $G_1^{-2}G_3^{-1}p_0$ to $G_1G_2^{-1} p_0$;
    \item $G_1G_2 p_0$ to $G_3p_0$;
    \item $G_1^{-1}G_2^{-1} p_0$ to $G_1G_2 p_0$;
    \item $G_1^{-2}G_2^{-1} p_0$ to $G_1^2G_2 p_0$.
  \end{enumerate}
\end{prop}

\begin{pf}
  We show a slightly stronger statement, namely in order to show that
  $G_2^{-1}g p_0 = h p_0$, we will exhibit $h^{-1} G_2^{-1} g$ as an
  explicit power of $G_1$.

  The result follows from the presentation of the group (strictly
  speaking, they only depend on the relations we know to hold, not on
  the fact that this really gives a presentation). For the sake of
  brevity, we use word notation.

  \begin{enumerate}
    \item $\bar3 1 \bar 2 3 = \bar2\bar1 2 1 \bar 2 \bar2 1 2 = \bar2\bar1 \bar2 \bar1\bar2 = 1$;
    \item $3\bar 1\bar2\bar1^3\bar3 = \bar2 12\cdot 121121\cdot \bar 1\cdot 121121\cdot 2 = 
      \bar2 (12121) (12121)1212 =\bar2^4 \bar 1 = \bar1$;
    \item $21\bar2\bar1 3 = 21\bar2\bar1\bar212 = \bar1$;
    \item $3\bar2\bar1 2 = Id$;
    \item $2\bar1\bar2\bar1^2\bar3 = 2(\bar1\bar2\bar1)^22 = 2(121)2 = \bar1$;
    \item $\bar3\bar212 = Id$;
    \item $\bar2\bar1\bar2\bar1\bar2 = 1$;
    \item $\bar2\bar1^2\bar2\bar1^2\bar2 = 1^2$.
  \end{enumerate}
\end{pf}

On the level of vertices, we have
\begin{itemize}
  \item $G_2^{-1} p_2 = p_2$;
  \item $G_2^{-1} p_{\bar121} = p_{\bar323}$;
  \item $G_2^{-1} p_{21^3} = p_{23^3}$;
  \item $G_2^{-1} p_{121^2} = p_{323^2} = p_{12\bar1}$.
\end{itemize}

\subsection{Faces paired by $G_3$}

The corresponding statement about the side-pairing map for the other
two base faces is the following.
\begin{prop}\label{prop:pairing3}
  The isometry $G_3^{-1}$ maps 
  \begin{enumerate}
    \item $G_2 p_0$ to $G_1^{-1}G_3^{-1} p_0$;
    \item $G_2^{-1} p_0$ to $G_2^{-1} p_0$;
    \item $G_1G_2 p_0$ to $G_1^2 G_2 p_0$;
    \item $G_1G_2^{-1} p_0$ to $G_1G_3^{-1}p_0$;
    \item $G_1G_3p_0$ to $G_1^3G_2 p_0$;
    \item $G_1^{-1}G_3 p_0$ to $G_1^{-1}G_2^{-1}p_0$.
  \end{enumerate}
\end{prop}

\begin{pf}
  \begin{enumerate}
    \item follows from $31\bar32=\bar2121\bar2\bar1\bar2\bar2=\bar2(\bar2\bar1\bar2)^2\bar2=1$;
    \item follows from $2\bar3\bar2=2\bar2\bar1 = 1$;
    \item follows from $\bar2\bar1^2\bar312 =  \bar2\bar1(212)^2=\bar2\bar1\bar2\bar1\bar2=1$;
    \item follows from $3\bar1\bar3\bar1\bar2 = Id$;
    \item follows from $\bar2\bar1^3\bar313=\bar2\bar1\bar1 \cdot\bar1\bar2\bar1\cdot 212 \cdot 22 12 = 
      \bar2\bar1(212)^2 = \bar2\bar1\bar2\bar1\bar2 = 1$;
    \item follows from $21\bar3 \bar1 3 = 21\bar2\bar12 \bar1\bar2 12 = 21 \bar2\bar1\bar2\bar2\bar2\bar1\bar2 12
      = (21212)^3 = \bar1^2$.
  \end{enumerate}
\end{pf}

On the level of vertices, we have
\begin{itemize}
  \item $G_3^{-1} p_2 = p_{\bar323}$;
  \item $G_3^{-1} p_{12\bar 1} = p_{1^32}$.
\end{itemize}
The last equality holds because
$$
\bar312\bar13 = \bar2\bar1212\bar1\bar212 = 1(212)^3112=1^32.
$$

\section{Ridge cycles} \label{sec:cycles}

Because of Giraud's theorem, the ridge cycles automatically satisfy
the hypotheses of the Poincar\'e polyhedron theorem. In particular, we
get the following:
\begin{thm}
  $D$ is a fundamental domain for the action of cosets of $\langle
  G_1\rangle$ in $\Gamma$. In particular, $D=F$ (see
  Theorem~\ref{thm:partialford}).
\end{thm}

Every ridge cycle is equivalent to one of the cycles listed in
Table~\ref{tab:cycles} (equivalent means that we allow shifting within
the cycle, and also conjugation by a power of $G_1$). We list the
cycle until we come back to the image of the initial ridge under a
power $G_1^k$ (in that case we close up the cycle by $G_1^{-k}$).
\begin{table}[htbp]
$$
\begin{array}{|c|c|}
  \hline
  2\cap 3 
      \stackrel{\bar 2}{\longrightarrow} \bar131\cap\bar2 
      \stackrel{\bar 1\bar31}{\longrightarrow} \bar3\cap\bar1\bar31 
      \stackrel{3}{\longrightarrow} 2\cap 3 
             & 2=[3,\bar1]\\   
  2\cap \bar1^3\bar31^3 
      \stackrel{\bar 2}{\longrightarrow} 1\bar3\bar1\cap\bar2 
      \stackrel{13\bar 1}{\longrightarrow} 3\cap 13\bar1 
      \stackrel{\bar 3}{\longrightarrow} 1^32\bar1^3\cap\bar3  
             & \bar1^3\bar313\bar1\bar2 = id\\ 
  2\cap \bar131 
      \stackrel{\bar 2}{\longrightarrow} \bar1\bar21\cap\bar2
      \stackrel{\bar 121}{\longrightarrow} \bar1^3\bar31^3\cap \bar121
      \stackrel{\bar 1^331^3}{\longrightarrow} \bar1^221^2\cap \bar1^331^3
             & \bar131^221\bar2 = id\\
  2\cap \bar121 \stackrel{\bar 2}{\longrightarrow} \bar3\cap \bar2
      \stackrel{3}{\longrightarrow} \bar2\cap 3
      \stackrel{2}{\longrightarrow} 12\bar1 \cap 2
             & 12=23\\
  2\cap \bar1\bar2 1 \stackrel{\bar 2}{\longrightarrow} 12\bar1\cap \bar2
             & (12)^3\\
  2\cap \bar1^2\bar21^2 \stackrel{\bar 2}{\longrightarrow} 1^22\bar1^2\cap \bar2
             & (121)^3\\\hline
\end{array}
$$
\caption{Ridge cycles, and the corresponding relation in the group.}\label{tab:cycles}
\end{table}

Using the relations
$$
12=23,\quad (12)^3=(121)^3=id,
$$
the other relations give $2^4=id$. Indeed, $\bar1^3\bar313\bar1\bar2=id$ gives
$$
id=\bar1^2 \bar313\bar1\bar2 \bar1 = \bar1^2 \bar2\bar12\cdot 1 \bar212\cdot \bar1\bar2\bar1 = \bar 1 (121)^2 2 1\bar212(121)^2 = 21\bar2^3121 = 21 (\bar2^4) \bar1\bar2.
$$
It is easy to check that the above set of relations is actually \emph{equivalent} to 
$$
12=23,\quad (12)^3=(121)^3=2^4=id.
$$

We summarize the above discussion in the following:
\begin{thm}\label{thm:derauxfalbel}
  The group $\Gamma$ has a presentation given by
  $$
  \langle G_1,G_2,G_3\ |\ G_2=[G_3,G_1^{-1}], G_1G_2=G_2G_3, G_2^4=id, (G_1G_2)3=id, (G_2G_1G_2)^3=id\rangle
  $$
\end{thm}

\section{Topology of the manifold at infinity} \label{sec:topology}

In this section, we prove that $\Gamma\setminus\Omega$ is indeed
homeomorphic to the figure eight knot complement. This was already
proved in~\cite{derauxfalbel} using a very different fundamental
domain for the action of the group.

We write $F$ for the Ford domain for $\Gamma$, $E$ for
$\partial_\infty F$, and $C$ for $\partial E$. By construction $F$,
$E$ and $C$ are all $G_1$-invariant.  

We will use Heisenberg coordinates $(z,t)$ for $\partial \ch
2\setminus \{p_\infty\}$, see section~\ref{sec:finite}. In these
coordinates, the action of $G_1$ is given by
\begin{equation}\label{eq:g1heis}
G_1(z,t)=(z-1,t+\Im(z)).
\end{equation}

It follows from the results in section~\ref{sec:statecombo} that $C$
is tiled by hexagons, and that there are four orbits of these hexagons
under the action of $G_1$. We need a bit more information about the
identifications on these hexagons, namely we need
\begin{itemize}
 \item The incidence relations between various hexagons, and
 \item The identifications on $C$ given by side-pairing maps.
\end{itemize}
The indicence relations follow immediately from the results in
section~\ref{sec:statecombo}, and it is summarized in
Figure~\ref{fig:hexagons}.

The union $U$ of the four hexagons labelled 1,2,3,4 is embedded in
$C$, and the action of $G_1$ induces identifications on $\partial
U$. We denote by $\sim$ the corresponding equivalence relation on $U$;
it is easy to check that $U/\sim$ is a torus.

We get the following result.
\begin{prop}
  $C$ is an unknotted topological cylinder, and $E$ is the region
  exterior to $C$.
\end{prop}

\begin{pf}
It follows from the fact that $C$ is invariant under the action of
$G_1$ that it is an unknotted cylinder in $\C\times\R$ (it is a
$\Z$-covering of $C/\langle G_1\rangle$). In fact, the real axis gives a core curve
for the solid cylinder bounded by $C$. In view of $G_1$-invariance, it
is enough to check that the interval $[0,1]$ on the $x$-axis is
outside $E$. This is readily checked, in fact this interval is
actually completely inside the spinal sphere $\Ss_1$.  
\end{pf}

The identifications in $C$ come from side pairings, which are
described in section~\ref{sec:pairings}. Figures~\ref{fig:face2}
and~\ref{fig:face3} contain a list of vertices, which are uniquely
determined by the list of faces they are on (in fact they are on
precisely three bisectors).

For instance, there is a vertex on $b_1\cap b_3\cap G_1(b_1)$. By
Proposition~\ref{prop:pairing2}, $G_2^{-1}$ maps this to the vertex on
$b_2\cap b_3\cap G_1^{-1}b_3$. The vertex on $b_1\cap b_3\cap
G_1^{-1}(b_3)$ is mapped to the vertex on $b_2\cap G_1^{-1}(b_2)\cap
G_1^{-1}b_3$. The image of these two points determine the image of the
entire hexagon on $b_1$ (in Figure~\ref{fig:hexagons}, the map flips
the orientation of the hexagon).

By doing similar verifications, one checks that the identification
pattern on the hexagons on $\Ss_1,\dots,\Ss_4$ is the same as the one
for the Ford domain of the holonomy of the real hyperbolic structure
on the figure eight knot complement, see Figure~\ref{fig:torusLink3}.

Now since the exterior of $C$ is homeomorphic to $C\times
[0,+\infty[$ (in a $G_1$-equivariant way), we get:

\begin{cor}
  $\Gamma\setminus E$ is homeomorphic to the figure eight knot complement.
\end{cor}

\section{Stability of the combinatorics} \label{sec:stabcombo}

The first remark is that distinct bounding bisectors for the Ford
domain for the unipotent solution are never cospinal, and as a
consequence the intersections
$\widehat{\gamma}_1\cap\widehat{\gamma}_2$ are uniquely determined by
the triple $p_0,\gamma_1 p_0,\gamma_2 p_0$. Of course, this property
will hold for all values of the twist parameter of $G_1$.

Now every point of an open 2-face is on precisely two bounding
bisectors, and that intersection is transverse. In other words, every
open 2-face will survive in small perturbations.

A similar remark holds for 1-faces, namely no 1-face of the Ford
domain for the boundary unipotent case is contained in a geodesic. In
fact every point on an open 1-face is on precisely three bounding
bisectors, and these intersect transversely as well.

The only issue is to analyze vertices. There is nothing to check for
the ideal vertices, since they are defined as the intersection of four
hypersurfaces (three bounding bisectors and the boundary of the ball)
that intersect transversely.

The finite vertices are on more than four bounding bisectors, but they
are also fixed by elliptic elements in the group. In fact, we already
justified that they stayed on the same bisectors for small
deformations, see section~\ref{sec:fixedpoints}, more specifically
Lemmas~\ref{lem:p2} and~\ref{lem:p1b21}. The transversality statement
of Proposition~\ref{prop:transversevertices} will remain true for
small perturbations as well.

This implies that the combinatorics stay stable in small deformations.

\section{Stability of the side pairing}

Let $F^{(0)}$ be the Ford domain for the boundary unipotent group, and
$F^{(t)}$ the one for the twist parabolic group corresponding to
parameter $t$.

The proof that $F^{(0)}$ has side-pairings relies on the determination of
the precise combinatorics, and also of the group relations. By the
previous section, the combinatorics are stable, and by
Proposition~\ref{prop:twistrelations}, the relations hold throughout
the deformation. The proof of Propositions~\ref{prop:pairing2}
and~\ref{prop:pairing3} then shows that $F^{(t)}$ has side-pairings, at
least for small values of $t$.

The verification that the Ford domain for the boundary unipotent group
satisfies the hypotheses of the Poincar\'e polyhedron theorem is given
in section~\ref{sec:cycles}. Since all intersections of bounding
bisectors are Giraud disks, the cycle condition is a direct
consequence of the existence of pairings. 

Let $\Gamma_t$ denote the image of $\rho_t$. We now get:
\begin{thm}
  There exists a $\delta>0$ such that whenever $|t|<\delta$,
  $\Gamma_t$ is discrete with non-empty domain of discontinuity, its
  manifold at infinity is homeomorphic to the figure eight knot
  complement, and it has the presentation
    $$
     \langle G_1,G_2,G_3\ |\ G_2=[G_3,G_1^{-1}], G_1G_2=G_2G_3, G_2^4=id, (G_1G_2)^3=id, (G_2G_1G_2)^3=id\rangle.
    $$
\end{thm}

\bibliographystyle{plain}
\bibliography{biblio}

\begin{thebibliography}{10}

\bibitem{agg}
S.~Anan'in, C.~H. Grossi, and N.~Gusevskii.
\newblock Complex hyperbolic structures on disc bundles over surfaces.
\newblock {\em Int. Math. Res. Not. IMRN}, 19:4295--4375, 2011.

\bibitem{semialgebraic}
S.~Basu, R.~Pollack, and M.-F. Roy.
\newblock {\em Algorithms in real algebraic geometry}, volume~10 of {\em
  Algorithms and Computation in Mathematics}.
\newblock Springer-Verlag, Berlin, second edition, 2006.

\bibitem{bfg}
N.~Bergeron, E.~Falbel, and A.~Guilloux.
\newblock Tetrahedra of flags, volume and homology of {SL}(3).
\newblock {\em Geom. Topol.}, 18:1911--1971, 2014.

\bibitem{cohen}
H.~Cohen.
\newblock {\em A course in computational algebraic number theory}, volume 138
  of {\em Graduate Texts in Mathematics}.
\newblock Springer-Verlag, Berlin, 1993.

\bibitem{deraux4445}
M.~Deraux.
\newblock Deforming the {$\mathbb{R}$}-{F}uchsian (4,4,4)-triangle group into a
  lattice.
\newblock {\em Topology}, 45:989--1020, 2006.

\bibitem{derauxcensusfkr}
M.~Deraux.
\newblock On spherical {CR} uniformization of 3-manifolds.
\newblock {\em Exp. Math.}, 24:335--370, 2015.

\bibitem{derauxfalbel}
M.~Deraux and E.~Falbel.
\newblock Complex hyperbolic geometry of the figure eight knot.
\newblock {\em Geom. Top.}, 19:237--293, 2015.

\bibitem{dpp2}
M.~Deraux, J.~R. Parker, and J.~Paupert.
\newblock New non-arithmetic complex hyperbolic lattices.
\newblock To appear in {I}nvent. {M}ath., arXiv:1401.0308.

\bibitem{falbelFigure8}
E.~Falbel.
\newblock A spherical {CR} structure on the complement of the figure eight knot
  with discrete holonomy.
\newblock {\em J. Differential Geom.}, 79(1):69--110, 2008.

\bibitem{fgkrt}
E.~Falbel, A.~Guilloux, P.-V. Koseleff, F.~Rouillier, and M.~Thistlethwaite.
\newblock Character varieties for {SL}(3,{C}) : the figure eight knot.
\newblock Preprint, arXiv:1412.4711.

\bibitem{fkr}
E.~Falbel, P.-V. Koseleff, and F.~Rouillier.
\newblock Representations of fundamental groups of $3$-manifolds into
  {PGL}$(3,\mathbb{C})$: exact computations in low complexity.
\newblock Preprint, arXiv:1307.6697.

\bibitem{giraud}
G.~Giraud.
\newblock Sur certaines fonctions automorphes de deux variables.
\newblock {\em Ann. Sci. {\'E}cole Norm. Sup.}, 38:43--164, 1921.

\bibitem{goldmantorus}
W.~M. Goldman.
\newblock Conformally flat manifolds with nilpotent holonomy and the
  uniformization problem for {$3$}-manifolds.
\newblock {\em Trans. Amer. Math. Soc.}, 278(2):573--583, 1983.

\bibitem{goldman}
W.~M. Goldman.
\newblock {\em Complex {H}yperbolic {G}eometry}.
\newblock Oxford Mathematical Monographs. Oxford University Press, 1999.

\bibitem{gkl}
W.~M. Goldman, M.~Kapovich, and B.~Leeb.
\newblock Complex hyperbolic manifolds homotopy equivalent to a {R}iemann
  surface.
\newblock {\em Comm. Anal. Geom.}, 9(1):61--95, 2001.

\bibitem{goldmanparkerdirichlet}
W.~M. Goldman and J.~Parker.
\newblock Dirichlet polyhedra for dihedral groups acting on complex hyperbolic
  space.
\newblock {\em J. of Geom. Analysis}, 2:517--554, 1992.

\bibitem{hmp}
M.~Heusener, V.~Munoz, and J.~Porti.
\newblock The {SL}(3,{C})-character variety of the figure eight knot.
\newblock Preprint, arXiv:1505.04451.

\bibitem{kamishimatsuboi}
Y.~Kamishima and T.~Tsuboi.
\newblock C{R}-structures on {S}eifert manifolds.
\newblock {\em Invent. Math.}, 104(1):149--163, 1991.

\bibitem{JRP-book}
J.~R. Parker.
\newblock {\em Complex {H}yperbolic {K}leinian {G}roups}.
\newblock Cambridge University Press, To appear.

\bibitem{parkerplatis}
J.~R. Parker and I.~D. Platis.
\newblock Open sets of maximal dimension in complex hyperbolic quasi-{F}uchsian
  space.
\newblock {\em J. Differential Geom.}, 73(2):319--350, 2006.

\bibitem{parkerwill}
J.~R. Parker and P.~Will.
\newblock A complex hyperbolic {R}iley slice.
\newblock Preprint, arXiv:1510:01505.

\bibitem{rileyquadratic}
R.~Riley.
\newblock A quadratic parabolic group.
\newblock {\em Math. Proc. Cambridge Philos. Soc.}, 77:281--288, 1975.

\bibitem{riley2bridge}
R.~Riley.
\newblock Nonabelian representations of {$2$}-bridge knot groups.
\newblock {\em Quart. J. Math. Oxford Ser. (2)}, 35(138):191--208, 1984.

\bibitem{rouillier}
F.~Rouillier.
\newblock Solving zero-dimensional polynomial systems through the rational
  univariate representation.
\newblock Technical report, INRIA, 1998.

\bibitem{richWLC}
R.~E. Schwartz.
\newblock Degenerating the complex hyperbolic ideal triangle groups.
\newblock {\em Acta Math.}, 186(1):105--154, 2001.

\bibitem{schwartzICM}
R.~E. Schwartz.
\newblock Complex hyperbolic triangle groups.
\newblock In {\em Proceedings of the {I}nternational {C}ongress of
  {M}athematicians, {V}ol. {II} ({B}eijing, 2002)}, pages 339--349, Beijing,
  2002. Higher Ed. Press.

\bibitem{rhochi}
R.~E. Schwartz.
\newblock Real hyperbolic on the outside, complex hyperbolic on the inside.
\newblock {\em Inv. Math.}, 151(2):221--295, 2003.

\bibitem{richBook}
R.~E. Schwartz.
\newblock {\em Spherical {CR} geometry and {D}ehn surgery}, volume 165 of {\em
  Annals of Mathematics Studies}.
\newblock Princeton University Press, 2007.

\bibitem{toledoJDG}
D.~Toledo.
\newblock Representations of surface groups in complex hyperbolic space.
\newblock {\em J. Diff. Geom.}, 29:125--133, 1989.

\bibitem{will}
P.~Will.
\newblock The punctured torus and {L}agrangian triangle groups in {${\rm
  PU}(2,1)$}.
\newblock {\em J. Reine Angew. Math.}, 602:95--121, 2007.

\end{thebibliography}

\end{document}